\documentclass[12pt]{amsart}

\usepackage[dvipsnames]{xcolor}
\usepackage[ruled]{algorithm2e}
\usepackage{amsbsy}
\usepackage{amsfonts}
\usepackage{amsmath}
\usepackage{amssymb}
\usepackage{amsthm}
\usepackage{anyfontsize}
\usepackage[title]{appendix}
\usepackage[font = footnotesize]{caption}
\usepackage{enumitem}
\usepackage{float}
\usepackage[T1]{fontenc}
\usepackage{graphicx}
\usepackage{hyperref}
\usepackage[noabbrev]{cleveref}
\usepackage{mathtools}
\usepackage[numbers,sort&compress]{natbib}
\usepackage{nccmath}
\usepackage{nicematrix}
\usepackage{pgfplots}
\usepackage[font=footnotesize]{subcaption}
\usepackage{tikz}
\usepackage{verbatim}
\usepackage{wavelet}
\usepackage{xparse}
\usepackage[top=0.6in,bottom=0.6in,left=0.8in,right=0.8in]{geometry}

\DeclareMathOperator* {\Span} {span}

\Crefname{thm}{Theorem}{Theorems}
\Crefname{cor}{Corollary}{Corollaries}
\Crefname{prop}{Proposition}{Propositions}
\Crefname{example}{Example}{Examples}
\Crefname{lemma}{Lemma}{Lemmas}
\Crefname{figure}{Figure}{Figures}
\Crefname{appendix}{Appendix}{Appendices}

\newcommand {\ii} {\mathbf{i}} 
\newcommand {\boldone} {\mbox{\usefont{U}{bbold}{m}{n}1}}

\newcommand{\cut}[1]{ \textcolor{red}{} } 
\renewcommand{\SS}{\mathcal{S}} 
\newcommand{\spt}{\mathbf{c}^*} 
\newcommand{\bpt}{\mathbf{b}^*} 
\newcommand{\LN}{\mathsf{\Pi}}
\newcommand{\LNC}{\mathsf{\Gamma}}

\DeclareDocumentCommand {\dd} {mo}
{ \IfValueTF{#2}
    { \left( \frac{\mathrm{d}} {\mathrm{d} #1} \right)^{\!#2} \!\! }
    { \frac{\mathrm{d}} {\mathrm{d} #1} }
}

\makeatletter
\DeclareDocumentCommand {\itemequation} {mooo}
{%
    \vspace{0.5em}
    \IfValueTF{#3}{\item[#3]}{\item} %
    \begingroup
        \refstepcounter{equation} %
        \IfValueTF{#2}{#2 \sbox2{\@eqnnum}}{} %
        \IfValueTF{#4}{\sbox4{#4}}{} %
        \sbox1{$\displaystyle #1 \m@th$} %
        \dimen@=.5\dimexpr\linewidth-\wd1\relax
        \ifcase
            \ifdim\wd4>\dimen@ \z@
            \else
                \ifdim\wd2>\dimen@ \z@
                \else \@ne
                \fi 
            \fi
            \@latex@warning{Equation is too large}%
        \fi
        \noindent   
        \rlap{\copy4}%
        \rlap{\hbox to \linewidth{\hfill\copy1\hfill}}%
        \hbox to \linewidth{\hfill\copy2}%
        \vspace{0.5em} %
    \endgroup
    \ignorespaces 
}
\makeatother  

\newtheorem{thm}{Theorem}[section]

\newtheorem{prop}[thm]{Proposition}
\newtheorem{lemma}[thm]{Lemma}

\newtheorem{example}[thm]{Example}

\numberwithin{equation}{section}

\allowdisplaybreaks

\hypersetup {
    colorlinks = true,
    linkcolor = black,
    filecolor = magenta,
    urlcolor = NavyBlue,
    citecolor = NavyBlue
}

\title[An Efficient FDM-Based PML Technique for Acoustic Scattering Problems]{An Efficient Finite Difference-Based PML Technique for Acoustic Scattering Problems}

\author{Bin Han and Jiwoon Sim}

\address{Department of Mathematical and Statistical Sciences, University of Alberta, Edmonton, Alberta, Canada T6G 2G1.}

\email{bhan@ualberta.ca, jiwoon2@ualberta.ca}



\begin{document}

	\maketitle
	
	\vspace{-2em}
	\begin{abstract}
		The acoustic scattering problem is modeled by the exterior Helmholtz equation, which is challenging to solve due to both the unboundedness of the domain and the high dispersion error, known as the pollution effect. We develop high-order compact finite difference methods (FDMs) in polar coordinates to numerically solve the problem with multiple arbitrarily shaped scatterers. The unbounded domain is effectively truncated and compressed via perfectly matched layers (PMLs), while the pollution effect is handled by the high order of our method and a novel pollution minimization technique. This technique is easy to implement, rigorously proven to be effective and shows superior performance in our numerous numerical results. The FDMs we propose in regular polar coordinates achieve fourth consistency order. Yet, combined with exponential stretching and mesh refinement, we can reach sixth consistency order by slightly enlarging the stencil at certain locations. Our numerical examples demonstrate that the proposed FDMs are effective and robust under various wavenumbers, PML layer thickness and shapes of scatterers.
	\end{abstract}

    \section{Introduction}
    \label{sec:intro}
    
    \subsection{Background}
    
    We consider the exterior Helmholtz equation:
    \begin{equation}
        \label{eq:PDE}
        \begin{cases}
            \Delta u + \kappa^2 u = f &\mbox{in } \R^2 \bs \overline{D}, \\
            u = g &\mbox{on } \partial D, \\
            \lim\limits_{r \to \infty} \sqrt{r} \left( \frac{\partial u}{\partial r} - \ii \kappa u \right) = 0,
        \end{cases}
    \end{equation}
    which describes the propagation of wave scattered by the scatterer $D$ in the frequency domain \cite{linton2001handbook,williams1999fourier}. Here $D$ is a bounded domain in $\R^2$ and $\kappa$ is a constant representing the wavenumber. The last condition in \eqref{eq:PDE} is called the Sommerfeld radiation condition, which ensures the well-posedness of the equation.
    
    One reason \cref{eq:PDE} is considered challenging is that it is defined in an unbounded domain. Common methods to deal with such unboundedness include boundary element method \cite{nedelec2001acoustic,colton2013integral}, absorbing boundary conditions (ABCs) \cite{engquist1977absorbing,bayliss1980radiation,tsynkov1998numerical} and perfectly matched/absorbing layers (PMLs or PALs) \cite{bermudez2007optimal,bermudez2008exact,wang2025novel,yang2021truly}. The boundary element method makes use of the boundary integral representation of the solution. It is flexible with the geometry of the scatterer, but it results in a dense linear system after discretization. This leads to high computational cost, especially when the wavenumber is large. The latter two methods solve the problem in a finite computational domain. To set up an ABC, one typically uses the Dirichlet-to-Neumann map at the truncated boundary. As pointed out in \cite{tsynkov1998numerical}, an exact ABC is nonlocal, involving integral equations or pseudo differential operators that introduce an extra layer of complexity for numerically solving the solution. On the other hand, local ABCs (e.g., \cite{bayliss1980radiation,engquist1977absorbing}) may suffer from loss of accuracy. A PML introduces an artificial layer filled with fictitious media that absorbs the wave and eliminates artificial reflection at the end of this layer. This was first introduced by B\'{e}renger in his studies of electromagnetic waves \cite{berenger1994perfectly,berenger1996three,berenger2002perfectly}. By imposing the Dirichlet boundary condition at the end of this layer, we obtain an equation that is much easier to solve, whose solution either coincides with the original one \cite{bermudez2008exact,yang2021truly}, or has an exponentially small error \cite{lassas1998existence}, depending on the type of the layer. In this article, we incorporate the PML technique to handle the unboundedness of the problem \eqref{eq:PDE}.
    
    The second challenge of \eqref{eq:PDE} is that, if we consider a large wavenumber $\kappa$, the numerical solution suffers from a high amount of dispersion error, known as the pollution effect. Although pollution-free FDMs are available in 1D \cite{wang2014pollution,han2021dirac}, the pollution effect can not be completely avoided in multiple dimensions \cite{babuska1997pollution} except in special circumstances \cite{han2021dirac,wang2017pollution}. One way to reduce the pollution effect is to use high order numerical schemes. In general, finite element methods (FEMs) are able to produce an $L^2$ error of $C \kappa^{p + 1} h^p$, where $p$ is the polynomial degree and $C$ is independent of $\kappa$, $h$ and $p$ \cite{du2015preasymptotic,babuska1997pollution,deraemaeker1999dispersion}. As for FDMs, the pollution effect is quantified via a dispersion analysis, where one measures the difference between true and numerical wavenumbers. To ensure that such a difference is bounded, one needs to impose the boundedness of $\kappa^3 h^2$ for second order FDMs \cite{chen2005adaptive,chen2011optimal} and of $\kappa^5 h^4$ for fourth order FDMs \cite{dastour2019fourth}. These results consistently suggest that high order numerical schemes impose a less stringent requirement for the mesh size $h$. When it comes to high order FDMs, one usually seeks FDMs with compact stencils, that is, the stencil consisting of a 3 by 3 Cartesian product of grid points. This simplifies the implementation of the FDM and makes the resulting linear system more sparse, hence saving computational resources. As reported in \cite{han2025convergent,feng2025symmetric}, a compact FDM for the variable Poisson equation can achieve at most sixth consistency order. For Helmholtz equations, compact sixth order FDMs are achieved by \cite{turkel2013compact,wu2018optimal,zhang2019sixth}. It is also natural to consider the Helmholtz equation in polar coordinates. As the resulting PDE is no longer isotropic in $r$ and $\theta$ directions (see \eqref{eq:PDE_r}), one may expect a lower consistency order. In \cite{britt2010compact} the authors constructed a compact fourth order FDM for the Helmholtz equation in polar coordinates.
    
    We wish to further save computational effort by solving \cref{eq:PDE} using a mesh size at a scale of $\kappa^{-1}$, which is necessary to resolve a wave with wavelength $2\pi \kappa^{-1}$. In this case, high order methods alone only reduce the pollution effect to a moderate level, and additional measures must be taken. As FDMs with a prescribed consistency order often include free parameters, one is able to minimize the difference between true and numerical wavenumbers. This is done by minimizing the dispersion error obtained by substituting the plane wave solution with a discrete set of incoming angles into the FDM (e.g., \cite{chen2011optimal,wu2014dispersion,dastour2019fourth,dastour2021generalized,wu2018optimal}). Their respective consistency orders are 2, 2, 4, 4, 6, with stencil sizes of up to 9, 9, 25, 25, and 9. Recently, Feng et al. \cite{feng2021sixth} proposed a compact FDM with up to sixth consistency order for the interior Helmholtz equation with interfaces. In this work, the pollution effect is reduced by minimizing the averaged local truncation error obtained from the plane wave solution of all directions. These approaches effectively reduce the pollution effect to a large extent.
    
    \subsection{Contribution}
    In summary, we develop high-order compact FDMs combined with PMLs in polar coordinates. We propose a novel pollution minimization technique to handle the pollution effect. One notable difference (which is also a simplification) compared to all existing methods is that we do not compute a family of stencil coefficients with a prescribed order in advance, and then look for the best free parameters that minimize the pollution effect. In our approach, the minimization process applies to the stencil coefficients themselves, regardless of the free parameters. We will rigorously show that the minimizer automatically forms a set of stencil coefficients with prescribed consistency order as long as a suitable stencil is chosen. In contrast to traditional methods (e.g. \cite{feng2025symmetric,han2025convergent}) where one needs to perform heavy symbolic computation to obtain the stencil coefficients with long expressions, our method consists merely of solving a small-scale minimization problem.
    
    The FDMs we propose in regular polar coordinates achieve fourth consistency order, which is in par with \cite{britt2010compact}. Yet, we will combine it with exponential stretching and mesh refinement. By slightly enlarging the stencil at certain locations, we are able to achieve sixth consistency order. The mesh refinement also aids us in dealing with multiple non-circular scatterers. The numerical examples demonstrate the superior performance of our pollution minimization method, and verify that the proposed FDMs are effective and robust under various wavenumbers, PML layer thickness and shapes of scatterers.
    
    \subsection{Organization and notations}
    
    We first give a brief introduction of the PML technique in \Cref{sec:PML}. \Cref{sec:FDM} is devoted to the construction of FDMs as well as our method to minimize the pollution effect. \Cref{sec:mesh} describes the exponential stretching, mesh refinement, and how we handle non-circular scatterers. We present several numerical examples from diverse perspectives in \Cref{sec:numerical} to show the effectiveness and robustness of our method. Some mathematical proofs are gathered in \Cref{app:math}.
    
    For readers' convenience, we collect some commonly used notations below:
    \begin{itemize}
    	\item $\N = \{1, 2, 3, \ldots\}$, $\N_0 = \N \cup \{0\}$.
		\item $\td(a)$ takes value $1$ when $a = 0$; $0$ otherwise.
    	\item For a vector $\bm{k} = (k_1, k_2) \in \N_0^2$, we denote $|\bm{k}| = k_1 + k_2$. If $v$ is a smooth function in variables $(x_1, x_2)$, then its partial derivatives are denoted by $\partial^{\bm{k}} v := \partial_{x_1}^{k_1} \partial_{x_2}^{k_2} v$.
		\item The letter $h$ with a subscript represents the mesh size in this coordinate. $h$ itself is an indicator of the mesh size, whose precise definition varies according to the context.
		\item $\T$ denotes the unit circle. The function spaces $L^2(\T)$ and $H^m (\T)$ are the usual Lebesgue space and Sobolev space of $2\pi$-periodic functions. The norms of $L^2(\T)$ and $H^m (\T)$ are averaged, e.g., $\|f\|_{L^2(\T)} := ( \frac{1}{2\pi} \int_{-\pi}^{\pi} |f(t)|^2 \mathrm{d}t )^{1/2}$.
		\item The big-O notation $\bo(h^M)$ stands for a function that is bounded by $C h^M$ as $h \to 0^+$. $f \sim h^M$ means $f h^{-M}$ is bounded above and below by a positive constant as $h \to 0^+$.
    \end{itemize}

    \section{The perfectly matched layer (PML)}
    \label{sec:PML}

    In this section, we briefly introduce the PML technique that is widely used to truncate the unbounded domain. As we are mainly interested in circular PMLs, we suppose the scatterer $D$ contains the origin and turn the PDE \eqref{eq:PDE} into polar coordinates:
    \begin{equation}
        \label{eq:PDE_r}
        \begin{cases}
            u_{rr} + \frac{1}{r} u_r + \frac{1}{r^2} u_{\theta \theta} + \kappa^2 u = f &\mbox{in } \R^2 \bs \overline{D}, \\
            u = g &\mbox{on } \partial D, \\
            \lim\limits_{r \to \infty} \sqrt{r} \left( u_r - \ii \kappa u \right) = 0.
        \end{cases}
    \end{equation}
    Here and afterward we slightly abuse the notations by identifying the functions $u$, $f$, $g$ and the domain $D$ in different coordinate systems.

	Suppose $D$, $\supp f$ and the domain of physical interest are contained within $\{r < r_*\}$. The PML layer can then be set up in the annulus $\{ r_* < r < r_M \}$ for some $r_M > r_*$. A complex coordinate transform $r \mapsto \rho(r) \in \C$ is introduced in the layer, satisfying $\rho(r_*) = r_*$, $\rho$~is smooth and $\re \rho'$, $\im \rho'$ are increasing in the layer. Typical choices of function $\rho$ include:
    \begin{itemize}
        \item[(1)] $\rho(r) = r + \ii \alpha \left( \dfrac{r - r_*}{r_M - r_*} \right)^n$, $n \in \N$ \cite{collino1998optimizing,chen2005adaptive,singer2004perfectly};
        \item[(2)] $\rho(r) = r + \ii \alpha \log \dfrac{r_M - r_*}{r_M - r}$ \cite{bermudez2007optimal,bermudez2008exact};
        \item[(3)] $\rho(r) = r_* + (\alpha_1 + \ii \alpha_2) \dfrac{r - r_*}{r_M - r}$ \cite{yang2021truly};
    \end{itemize}
    where $\alpha$, $\alpha_1$, $\alpha_2$ are positive constants. We do not follow the convention of scaling $\rho$ by $\kappa^{-1}$. This issue will be further discussed in \Cref{sec:numerical}. We further extend the function $\rho$ to the interval $(0, r_M)$ such that $\rho(r) = r$ for $0 < r \leq r_*$. To conclude, we impose the following assumptions on $\rho$:
    \begin{itemize}
        \item[($A_\rho$)] $\rho(r) = r$ for $0 < r \leq r_*$, $\rho |_{[r_*, r_M)}$ is real analytic, and $\re \rho'$, $\im \rho' > 0$ in $[r_*, r_M)$.
    \end{itemize}
    The analyticity of $\rho$ is required for some technical reasons in \Cref{thm:pollution}. Note that the extended $\rho$ may be nonsmooth at $r_*$. 

    The PML technique aims to approximate $u(\rho(r), \theta)$ by another function $v(r, \theta)$ defined on $\{ 0 < r < r_M \} \bs \overline{D}$. As it is demonstrated in \Cref{sec:estimates}, the solution $u(\rho(r), \theta)$ together with its derivatives decay exponentially in the PML layer. Hence, it is appropriate to add the homogeneous Dirichlet boundary condition at the outermost boundary $r = r_M$ for the function $v$. More precisely, we consider the following PDE for $v$:
    \begin{equation}
        \label{eq:PDE_PML_r}
        \begin{cases}
            v_{rr} + \frac{(\rho')^2}{\rho^2} v_{\theta \theta} + \left( \frac{\rho'}{\rho} - \frac{\rho''}{\rho'} \right) v_r + \kappa^2 (\rho')^2 v = (\rho')^2 f &\mbox{in } \Omega \cup \Omega^{\mathrm{PML}}, \\
            v_r^+ = \rho'(r_*) v_r^- &\mbox{on } \Gamma, \\
            v = g &\mbox{on } \partial D, \\
            v = 0 &\mbox{on } \{ r = r_M \}.
        \end{cases}
    \end{equation}
    Here $v_r^\pm (r^*, \theta) = \lim_{r \to r_*^\pm} v_r (r, \theta)$ and we denote the regular region $\Omega := \{ r < r_* \} \bs \overline{D}$, PML region $\Omega^{\mathsf{PML}} := \{ r_* < r < r_M \}$ and the interface $\Gamma := \{ r = r_* \}$. Note that $v$ is continuous across the interface. 
    
    As we can see, solving \cref{eq:PDE} via PML introduces an error with two parts: the difference between $u$ and $v$, and the error of numerically solving for $v$. It is shown in \cite{lassas1998existence} that, if $\im \rho_M < \infty$, then $v$ approximates $u$ exponentially in $\Omega$ as the thickness of the PML layer increases, and \cite{bermudez2008exact, yang2021truly} show that $u = v$ in $\Omega$ when $\im \rho_M = \infty$. On the other hand, letting $\im \rho_M = \infty$ will render the PDE \eqref{eq:PDE_PML_r} singular (e.g., one can prove $\limsup_{r \to r_M^-} \left| \tfrac{\rho'(r)}{\rho(r)} \right| = \limsup_{r \to r_M^-} \left| \tfrac{\rho''(r)}{\rho'(r)} \right| = \infty$). Readers can refer to \cite{collino1998optimizing,modave2014optimizing} to see how to balance these two parts of errors by optimizing the complex coordinate transform $\rho$. To avoid complexity, we choose $\rho$ such that $\im \rho_M < \infty$. The detailed choice of $\rho$ for reducing the first error is discussed in \Cref{sec:numerical}. For the rest of this article, we are mainly interested in developing FDMs to reduce the second error. 

    \section{Derivation of finite difference methods}
    \label{sec:FDM}
    
    In the following, we will develop FDMs to solve the Helmholtz-PML equation \eqref{eq:PDE_PML_r}. First, in \Cref{sec:FDM_generic} we review an existing but powerful method to develop a generic FDM for our problem. The next section introduces our technique to reduce the pollution effect, which automatically maintains the consistency order and is extremely simple to implement. For now we only consider a circular scatterer so that we do not need to care about the boundary. The treatment of non-circular scatterers needs to be postponed until \Cref{sec:non_circular}.
    
    In view of \cref{eq:PDE_PML_r} (also \eqref{eq:PDE_PML_s} later on), we study FDMs for PDEs in the following form:
    \begin{equation}
        \label{eq:PDE_PML}
        \begin{cases}
            v_{rr} = a(r) v_{\theta \theta} + b(r) v_r + \tilde{\kappa}(r) v + \tilde{f} &\mbox{in } \Omega \cup \Omega^{\mathrm{PML}}, \\
            v_r^+ = \beta v_r^- &\mbox{on } \Gamma, \\
            v = g &\mbox{on } \partial D, \\
            v = 0 &\mbox{on } \{ r = r_M \},
        \end{cases}
    \end{equation}
    where the constant $\beta \in \C$ replaces $\rho'(r_*)$. Due to assumption ($A_\rho$), the functions $a$, $b$ and $\tilde{\kappa}$ is smooth in $(0, r_*)$ and $(r_*, r_M)$, and at $r_*$ we have $a'(r_*^+) = \beta a'(r_*^-)$ (similar for $b$ and $\tilde{\kappa}$).
    
    We put a rectangular grid on the computational domain. As we are free to set the location of the PML layer, we can simply let the grid points fall on $\Gamma$ and $\{ r = r_M \}$. When we consider a circular scatterer $D$, we also suppose $\partial D$ coincides with one of the grid lines. In this way, the treatment at the interface is simpler, and we do not need to separately consider boundary stencils near $\{ r = r_M \}$ and $\partial D$. The mesh size in $r$ and $\theta$ direction are denoted as $h_r > 0$ and $h_\theta > 0$. We also put $h = h_r = \gamma h_\theta$. The set of nodes is given by $\Omega_h := (h_r \Z) \times (h_\theta \Z) \cap (\Omega \cup \Gamma \cup \Omega^{\mathsf{PML}})$. 
    
    \subsection{The generic method}
    \label{sec:FDM_generic}
    
    As we pursue compact stencils, that is, 9-point stencils consisting of the stencil center and 8 other points surrounding it, there are two possibilities: either the stencil lies within $\overline{\Omega}$, $\overline{\Omega^{\mathsf{PML}}}$, or the stencil center $\spt$ is on the interface $\Gamma$. When the stencil is entirely included in the regular or PML region, the FDM exactly follows the generic method in \cite{feng2025symmetric} since the jump condition (second identity in \eqref{eq:PDE_PML}) is not relevant in this situation. We just list the essential steps and differences without details. 
    
    As in \cite{feng2025symmetric}, we can obtain the following expansion of the solution $v$:
    \begin{equation}
        \label{eq:v_taylor_int}
        v(\spt + ph)
        = \sum_{\bm{\ell} \in \LN_{M + 1}} \sum_{k = |\bm{\ell}|}^{M + 1}
        A^k_{\bm{\ell}} (p) h^k \cdot \partial^{\bm{\ell}} v(\spt) + F(p) + \bo (h^{M + 2}), \quad \spt \in \Omega_h.
    \end{equation}
    Here $M \in \N$ is the consistency order we would like to achieve, $\LN_{M + 1} := \{ \bm{k} = (k_1, k_2) \in \N_0^2: |\bm{k}| \leq M + 1, k_1 \leq 1 \}$, and $A^k_{\bm{\ell}} (p)$ and $F(p)$ are explicitly known quantities involving derivatives of $a$, $b$, $\tilde{\kappa}$ and $\tilde{f}$ at $\spt$. The only difference is that
    \begin{equation*}
        \tilde{a}^{\bm{k}}_{\bm{\ell}} = 
        \begin{cases}
            \displaystyle
            \binom{k_1 - 2}{\ell_1} \partial^{k_1 - \ell_1 - 2} a, & \mbox{if} \ \ell_2 = k_2 + 2, \\
            \displaystyle
            \binom{k_1 - 2}{\ell_1 - 1} \partial^{k_1 - \ell_1 - 1} b + \binom{k_1 - 2}{\ell_1} \partial^{k_1 - \ell_1 - 2} \tilde{\kappa}, & \mbox{if} \ \ell_2 = k_2.
        \end{cases}
    \end{equation*}
	instead of \cite[equation (6.3)]{feng2025symmetric}. Here we take the convention that $\binom{n}{m} = 0$ for $n, m \in \Z$ such that $m < 0$ or $m > n$.
	
	When $\spt \in \Gamma$, writing $p = (p_1, p_2) \in \R^2$, we need to change $\partial^{\bm{\ell}} v(\spt)$ into $\partial^{\bm{\ell}} v^{\pm} (\spt)$ when $p_1 \neq 0$ and $\mathrm{sign}(p_1) = \pm 1$. In the case of $p_1 = 0$, for one thing, \cref{eq:v_taylor_int} can be equivalently derived from performing Taylor expansion solely in $\theta$ direction. In this way, we obtain $A^k_{\bm{\ell}} (p) = 0$ when $p_1 = 0$ and $\ell_1 = 1$, where $\bm{\ell} = (\ell_1, \ell_2)$. For another, it is clear that $\partial^{\bm{\ell}} v^+(\spt) = \partial^{\bm{\ell}} v^-(\spt)$ when $\ell_1 = 0$. Hence, it does not matter whether we choose $\partial^{\bm{\ell}} v^+ (\spt)$ or $\partial^{\bm{\ell}} v^- (\spt)$ in \eqref{eq:v_taylor_int} for $p_1 = 0$. Now, we use the jump condition in \eqref{eq:PDE_PML} to obtain $\partial^{\bm{\ell}} v^+(\spt) = \beta \partial^{\bm{\ell}} v^-(\spt)$ for all $\bm{\ell} = (\ell_1, \ell_2)$ with $\ell_1 = 1$. In conclusion, \eqref{eq:v_taylor_int} still holds when $\spt \in \Gamma$ by replacing $\partial^{\bm{\ell}} v(\spt)$ with $\partial^{\bm{\ell}} v^-(\spt)$ and using different expressions for $A^k_{\bm{\ell}} (p)$.
	
	The FDM we seek is in the form of
	\begin{equation}
		\label{eq:FDM}
		\mathcal{L}_h v_h (\spt) = \sum_{p \in \SS} C_p (\spt) F(p)
	\end{equation}
	where $\mathcal{L}_h v := \sum_{p \in \SS} C_p(\spt) v(\spt + ph)$ is the discretization operator, $\SS = \{ -1, 0, 1 \} \times \{ -\gamma^{-1}, 0, \gamma^{-1} \}$ is the reference stencil and $C_p (\spt) = \sum_{j = 0}^{M + 1} c_{p, j} h^j$ is the stencil coefficient on $\spt + ph$. According to \cite[Lemma 6.1]{feng2025symmetric}, the coefficients $c_{p, j}$ are solved from the linear system
	\begin{equation}
		\label{eq:c_p,j}
		\sum_{p \in \SS} A^{|\bm{\ell}|}_{\bm{\ell}}(p) c_{p, j}
		= -\sum_{k = 0}^{j - 1} \sum_{p \in \SS} A^{|\bm{\ell}| + j - k}_{\bm{\ell}}(p) c_{p, k}, \quad \forall \, \bm{\ell} \in \LN_{M + 1 - j},
	\end{equation}
	iteratively for $j = 0, \ldots, M + 1$. Such stencil coefficients (if exist) lead to an $M$-th order consistent FDM, indicated below: 
	\begin{equation}
		\label{eq:FDM_error}
		\mathcal{L}_h v (\spt) = \sum_{p \in \SS} C_p (\spt) F(p) + \bo(h^{M + 2})
		\quad \text{for all} \ v \ \text{satisfying} \ \eqref{eq:PDE_PML}.
	\end{equation}
	A mild non-degenerate condition is imposed on $c_{p, 0}$ according to \cite{feng2025symmetric} to truly reach $M$-th consistency order. For convenience, one can enforce a necessary condition that $c_{p, 0}$ is not identically $0$ for $p \in \SS$.
	
	According to our symbolic computation, the linear system \eqref{eq:c_p,j} with $j = 0$ has a nontrivial solution for $M \leq 3$. This upper limit can be improved to $4$ when $\beta = 1$, i.e., $\rho \in C^1(0, r_M)$, but can not be improved further even if we set $h_r = h_\theta$. The fourth consistency order is in accord with existing results \cite{britt2010compact} on the compact fourth order FDM for the Helmholtz equation in polar coordinates. 

    \subsection{Pollution minimizing technique}
    \label{sec:pollution}
    
    Now we aim to resolve the pollution effect by minimizing the truncation error $\bo(h^{M + 2})$ in \eqref{eq:FDM_error} suitably. For any function $w$ such that $\Delta w + \kappa^2 w = 0$, the function $\tilde{f}$ in \eqref{eq:PDE_PML} is zero. Since the term $F(p)$ in \eqref{eq:v_taylor_int} is a linear combination of $\partial^{\bm{k}} \tilde{f}$ for $|\bm{k}| < M$ (see \cite{feng2025symmetric} for more details), we have $F(p) = 0$. If $\mathcal{L}_h$ is an $M$-th order consistent FDM in the sense of \eqref{eq:FDM_error}, then $\mathcal{L}_h w(\spt) = \bo(h^{M + 2})$. Therefore, we aim to determine the stencil coefficients through the following minimization problem:
    \begin{equation}
    	\label{eq:hat_C_p}
    	(\widehat{C}_p(\spt))_{p \in \SS} := \arg\min \{ \mathcal{I}_h(\vec{a}): \vec{a} = (a_p)_{p \in \SS} \in \R^{\# \SS}, a_{(0, 0)} = 1 \},
    \end{equation}
    where $\mathcal{I}_h(\vec{a})$ is related to $\big| \hspace{-2pt} \sum_{p \in \SS} a_p w(\spt + ph) \big|$ (defined later in \eqref{eq:I_h}). In \eqref{eq:hat_C_p} we have adopted a normalization condition $a_{(0, 0)} = 1$, since multiplying a constant to all stencil coefficients does not affect the accuracy of the FDM but changes the value of $I(\vec{a})$.
    
	The common approach for the test function $w$ is the plane incident wave $w(r, \theta; \theta_0) := \exp(\ii \kappa r\cos$ $(\theta - \theta_0))$ coming in an angle of $\theta_0$. However, this is no longer appropriate in the context of PMLs, as it grows out of control due to the complex coordinate transform $\rho(r)$. To this end, we redefine the function $w$ in the following way:
    \begin{equation}
    	\label{eq:pollution_w}
    	w(r, \theta; \theta_0) := 
    	\begin{cases}
    		\exp(\ii \kappa \rho(r) \cos(\theta - \theta_0)) &\mbox{in } \Omega, \\
            \displaystyle \sum_{j \in \Z} \ii^j \frac{J_j(\kappa r_*)}{H^{(1)}_j (\kappa r_*)} H^{(1)}_j (\kappa \rho(r)) e^{\ii j (\theta - \theta_0)} &\mbox{in } \Omega^{\mathsf{PML}},
    	\end{cases}
    \end{equation}
    where $J_j$ is the Bessel function and $H_j^{(1)}$ is the Hankel function of the first kind. One can refer to \cite{olver2010nist} for the definition and properties of these functions. The function $w$ in \eqref{eq:pollution_w} is well defined and satisfies the first equation in \eqref{eq:PDE_PML_r}, which can be easily proved following \cite[Theorem 7.2]{bermudez2008exact}. We follow \cite{feng2021Helmholtz} for the objective function $\mathcal{I}_h$ with an additional penalty term, that is,
    \begin{equation}
    	\label{eq:I_h}
    	\widetilde{\mathcal{I}}_h (\vec{a}) = \bigg\| \sum_{p \in \SS} a_p w(\spt + ph; \cdot) \bigg\|_{L^2 (\T)}^2, \quad
    	\mathcal{I}_h(\vec{a}) := \widetilde{\mathcal{I}}_h (\vec{a}) + \delta_h |\vec{a}|_{\ell^2}^2.
    \end{equation}
    For the moment, we can think of $\delta = \bo(h^{2M + 4})$ such that the penalty term preserves the accuracy order and stabilizes the results. Further discussion can be seen in \Cref{thm:pollution}, \Cref{sec:choice_param} and \Cref{app:float_error}.
        
	Let us explain our choice of $w$ in detail. First, let $u$ be the solution to the original Helmholtz equation \eqref{eq:PDE}. By \cite[Theorem 7.2]{bermudez2008exact}, $u$ can be represented as a series in $\Omega^{\mathsf{PML}}$:
	\begin{equation}
		\label{eq:u_series}
		u(\rho(r), \theta) = \sum_{j \in \Z} \frac{a_j}{H^{(1)}_j (\kappa r_*)} H^{(1)}_j (\kappa \rho(r)) e^{\ii j \theta},
	\end{equation}
	where $a_j = \langle e^{\ii j \cdot}, u|_\Gamma \rangle_{L^2 (\T)}$ is the Fourier coefficient of $u|_\Gamma$. This shows that $w$ and $u$ have the same form in $\Omega^{\mathsf{PML}}$. We also expect that they have the same decaying property in view of \Cref{prop:sol_decay}. Besides, $w$ has a consistent format in $\Omega$ and $\Omega^{\mathsf{PML}}$ due to
	\begin{equation}
		\label{eq:plane_wave_expansion}
		\exp(\ii \kappa \rho(r) \cos(\theta - \theta_0))
		= \sum_{j \in \Z} \ii^j \frac{J_j(\kappa r_*)}{J_j(\kappa r_*)} J_j (\kappa \rho(r)) e^{\ii j (\theta - \theta_0)}.
	\end{equation}
	This identity can be proved by regarding the right hand side of \eqref{eq:plane_wave_expansion} as a Fourier series and applying \cite[equation (10.9.2)]{olver2010nist}.
	
	Since $\mathcal{I}_h$ in \eqref{eq:I_h} is a nonnegative quadratic form in $\vec{a}$, the minimization problem \eqref{eq:hat_C_p} is equivalent to solving a linear system. Besides, the evaluation of integrals in \eqref{eq:I_h} can be replaced by a truncated summation in the Fourier domain, which converges much faster due to the super-exponential decay of Bessel functions $J_j$ with respect to its order $j$ \cite[equation (10.19.1)]{olver2010nist}. We summarize the procedure in the algorithm below. Due to the uniformity in $\theta$ direction, we only need to perform this algorithm at each different $r$ value, instead of each stencil point. 
	
	\SetKw{KwSolve}{solve}
	\SetKw{KwFrom}{from}
	\begin{algorithm}[t]
		\label{alg:pollution}
		\DontPrintSemicolon
		\caption{Stencil coefficients with reduced pollution effect}
		\KwIn{$\spt = (r, \theta)$, $\SS = \{ -1, 0, 1 \} \times \{ -\gamma^{-1}, 0, \gamma^{-1} \}$, $\rho$, $h$, $\delta_h$, threshold $J \in \N$}
		\KwOut{Stencil coefficients $\widehat{C}_p (\spt)$ satisfying \cref{eq:hat_C_p,eq:pollution_w,eq:I_h}}
		
		\For{$p = (p_1, p_2) \in \SS$, $|j| \leq J$}{
			\eIf{$r \leq r_*$}{
				$\gamma_{p, j} \gets \dfrac{J_j (\kappa r_*)}{J_j (\kappa r_*)} J_j (\kappa \rho(r + p_1h)) e^{\ii j p_2 h}$
			}{
				$\gamma_{p, j} \gets \dfrac{J_j (\kappa r_*)}{H^{(1)}_j (\kappa r_*)} H^{(1)}_j (\kappa \rho(r + p_1h)) e^{\ii j p_2 h}$
			}
		}
		\For{$p, q \in \SS$}{
			$w_{p, q} \gets \sum_{|j| \leq J} \overline{\gamma_{p, j}} \gamma_{q, j} + \delta_h \td(p - q)$
		}
		\KwSolve{$\vec{a}$} \KwFrom{$(w_{p, q})_{p, q \in \SS \bs \{(0, 0)\}} \, \vec{a} = -(w_{p, (0, 0)})_{p \in \SS \bs \{(0, 0)\}}$}\;
		$\widehat{C}_{(0, 0)} (\spt) \gets 1$, $(\widehat{C}_p (\spt))_{p \in \SS \bs \{(0, 0)\}} \gets \vec{a}$\;
	\end{algorithm}
	
	Although we only minimize the truncation error with a particular family of solutions $w(\cdot; \theta_0)$ to the Helmholtz equation, the following theorem shows that the consistency order is preserved for all solutions to the PDE. The proof is provided in \Cref{sec:proof}. 
	\newpage
	\begin{thm}
		\label{thm:pollution}
		Let $\SS \subseteq \R^2$, $M \in \N$ and $F(p)$ be defined in \eqref{eq:v_taylor_int}. Suppose assumptions (i), (ii), (v$_1$), or assumptions (i)-(iv) and (v$_2$) from below are satisfied.
		\begin{enumerate}[leftmargin = 2em]
			\item[(i)] The function $\rho$ satisfies assumption ($A_\rho$).
			
			\item[(ii)] There exists an $M$th order consistent FDM $\mathcal{L}_h v_h = \sum_{p \in \SS} C_p(\spt) F(p)$ with $\mathcal{L}_h v := \sum_{p \in \SS} C_p(\spt)$ $v(\spt + ph)$, i.e., \cref{eq:FDM_error} holds.
						
			\item[(iii)] $\displaystyle \min_{|\vec{a}|_{\ell^2} = 1} \widetilde{I}_h (\vec{a}) > 0$ for all $h$ in a neighbourhood of $0$, where $\widetilde{I}_h (\vec{a})$ is defined in \eqref{eq:I_h}.
						
			\item[(iv)] The linear system \eqref{eq:c_p,j} with $j = 0$ has a unique solution $(c_{p, 0})_{p \in \SS}$ with $c_{(0, 0), 0} \neq 0$ up to scaling.
						
			\item[(v$_1$)] $\delta_h$ is analytic with respect to $h$ and $0 < \delta_h \sim h^{2M + 4}$.
			
			\item[(v$_2$)] $\delta_h$ is analytic with respect to $h$ and $0 \leq \delta_h = \bo(h^{2M + 4})$.
		\end{enumerate}
		Define a new set of stencil coefficients $(\widehat{C}_p(\spt))_{p \in \SS}$ through \cref{eq:hat_C_p,eq:pollution_w,eq:I_h} and let $\widehat{\mathcal{L}}_h$ be the corresponding discretization operator. Then the new FDM $\widehat{\mathcal{L}}_h v_h = \sum_{p \in \SS} \widehat{C}_p (\spt) F(p)$ is $M$-th order consistent in the sense that \eqref{eq:FDM_error} holds for $\widehat{\mathcal{L}}_h$ and $(\widehat{C}_p(\spt))_{p \in \SS}$.
	\end{thm}
	
	For assumption (iii) in \cref{thm:pollution}, we believe that $\min_{|\vec{a}|_{\ell^2} = 1} \widetilde{I}_h (\vec{a}) \sim h^{2M + 4}$ with $M = 4$ and this is observed numerically. For assumption (iv), we perform a symbolic computation to the zeroth order stencil coefficients $(c_{p, 0})_{p \in \SS}$ under the condition that $\rho \in C^1(0, r_M)$ and $\SS = \{-1, 0, 1\} \times \{-\gamma, 0, \gamma\}$. It turns out that $(c_{p, 0})_{p \in \SS}$ is unique up to a multiple of
	\begin{equation*}
	\begin{matrix}
		-\rho(r)^{-2} - \gamma^2 \rho'(r)^{-2} && -10 \rho(r)^{-2} + 2 \gamma^2 \rho'(r)^{-2} && -\rho(r)^{-2} - \gamma^2 \rho'(r)^{-2} \\
		2 \rho(r)^{-2} - 10 \gamma^2 \rho'(r)^{-2} && 20 \rho(r)^{-2} + 20 \gamma^2 \rho'(r)^{-2} && 2 \rho(r)^{-2} + 10 \gamma^2 \rho'(r)^{-2}. \\
		-\rho(r)^{-2} - \gamma^2 \rho'(r)^{-2} && -10 \rho(r)^{-2} + 2 \gamma^2 \rho'(r)^{-2} && -\rho(r)^{-2} - \gamma^2 \rho'(r)^{-2}
	\end{matrix}
	\end{equation*}
	Here, the position of the value $c_{p, 0}$ is the same as the position of $p \in \SS$. We will prove in \Cref{lem:center_nz} that under existing assumptions on $\rho$ we must have $c_{(0, 0), 0} \neq 0$. Hence, by following Algorithm 1 with $\delta_h = 0$, we obtain a fourth-order FDM with reduced pollution effect. 
	
	Note that the fourth order consistency is theoretical, and the floating point error will affect its performance due to the ill-conditioning of the minimization problem \eqref{eq:hat_C_p}. This is discussed in \Cref{app:float_error}.
	
    \section{Exponential stretching and mesh refinement}
    \label{sec:mesh}

    In this section, we will talk about several benefits of adopting an additional exponential stretching before the complex coordinate transform $\rho$. Together with suitable mesh refinement, we are able to produce a quasi-uniform mesh seen from the Cartesian coordinates. Such uniformity contributes to an even distribution of numerical error within the computational domain.

    \subsection{Description and benefits}
    \label{sec:mesh_descr}
    
    Let $r = e^s$, $r_* = e^{s_*}$, $r_M = e^{s_M}$, and the complex coordinate transform $\rho$ is performed after it. Then the corresponding PDEs of \eqref{eq:PDE_r} and \eqref{eq:PDE_PML_r} in $(s, \theta)$ coordinates become
    \begin{equation}
        \label{eq:PDE_s}
        \begin{cases}
            \Delta u + \kappa^2 e^{2s} u = e^{2s} f &\mbox{in } \R^2 \bs \overline{D}, \\
            u = g &\mbox{on } \partial D, \\
            \lim\limits_{r \to \infty} \sqrt{r} \left( \frac{\partial u}{\partial r} - \ii \kappa u \right) = 0
        \end{cases}
    \end{equation}
    and
    \begin{equation}
        \label{eq:PDE_PML_s}
        \begin{cases}
            v_{ss} + (\rho')^2 v_{\theta \theta} - \frac{\rho''}{\rho'} v_r + \kappa^2 (\rho')^2 e^{2 \rho} v = (\rho')^2 e^{2 \rho} f &\mbox{in } \Omega \cup \Omega^{\mathrm{PML}}, \\
            v_s^+ = \rho'(s_*) v_s^- &\mbox{on } \Gamma, \\
            v = g &\mbox{on } \partial D, \\
            v = 0 &\mbox{on } \{ s = s_M \}.
        \end{cases}
    \end{equation}
    respectively. Note that the Laplace operator reappeared in \eqref{eq:PDE_s} as well as in \eqref{eq:PDE_PML_s} when $s < s_*$. This enables us to use sixth order compact FDMs. In \cite{feng2025symmetric} it is proved that 6 is the highest consistency order for two dimensional variable Poisson equation. We expect that the same is true for our case, as the argument in \cite{feng2025symmetric} depends essentially on the leading term of the PDE.

\begin{figure}[t]
        \centering
		\includegraphics[width = 180pt]{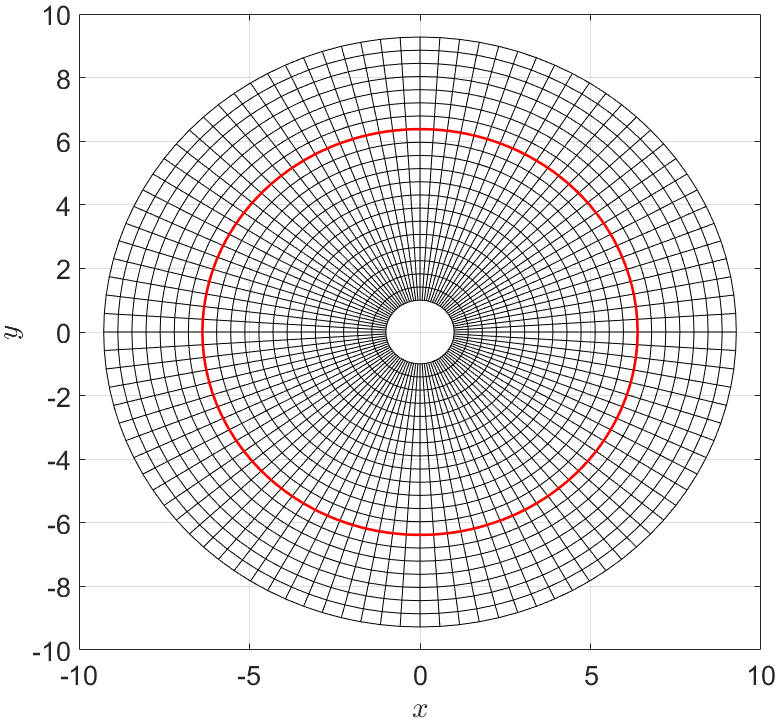}
		\includegraphics[width = 180pt]{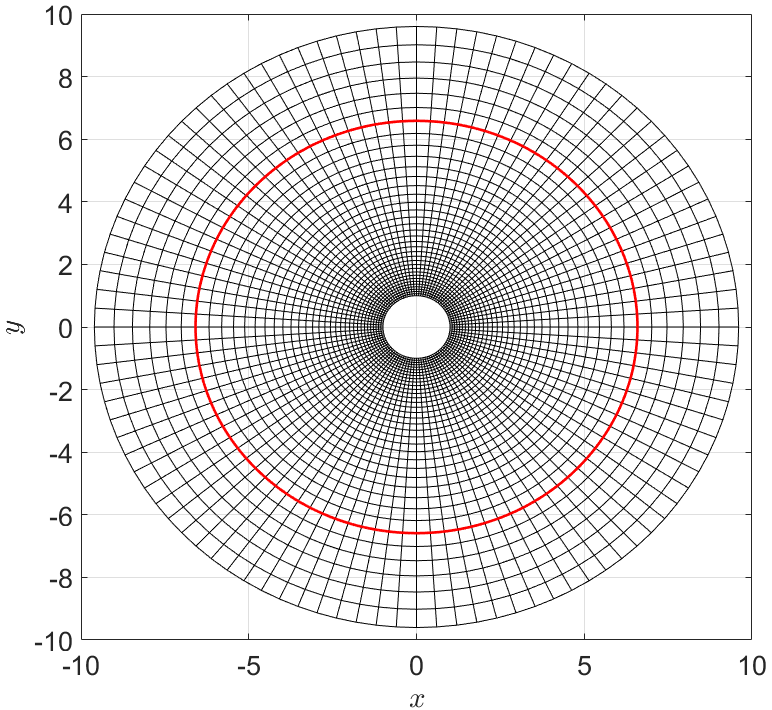}
		
		\includegraphics[width = 180pt]{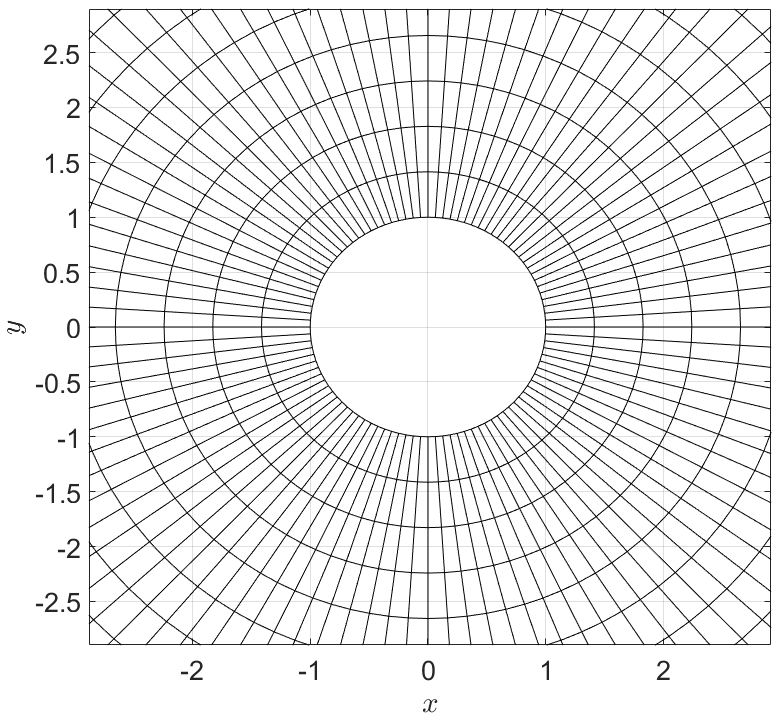}
		\includegraphics[width = 180pt]{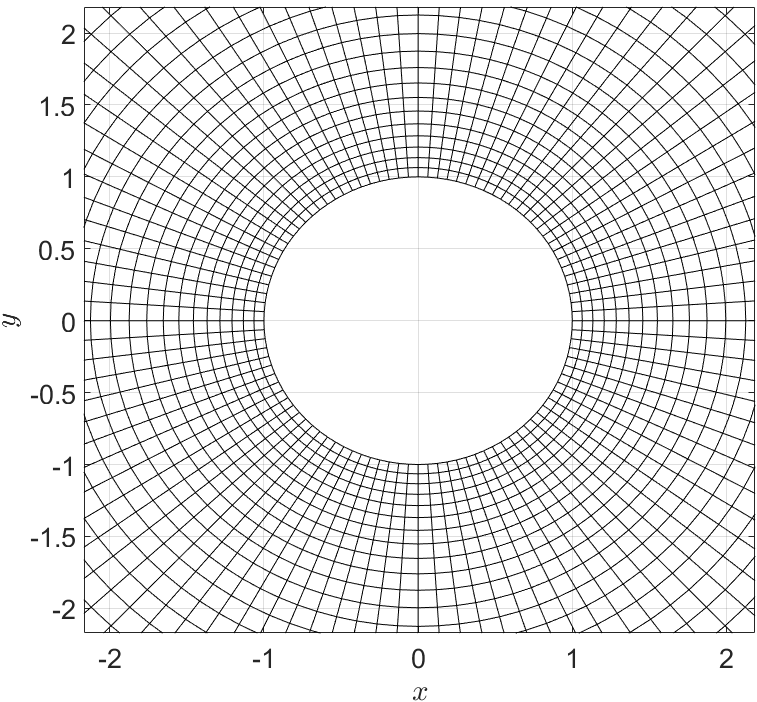}
        \caption{Comparison of meshes generated in $(r, \theta)$ (left) and $(s, \theta)$ coordinates (right). The mesh size at the interface is the same. The bottom figures are zoomed-in views of the corresponding top figures at the origin. Red circle indicates the interface $\Gamma$.}
        \label{fig:mesh}
    \end{figure}

    More importantly, let us compare the meshes generated in $(r, \theta)$ and $(s, \theta)$ coordinates. We put a rectangular mesh of size $h_r$, $h_\theta$ in the former coordinates and put a square mesh (required for sixth order FDM) of size $h_s$ for the latter. In the meantime, we keep the size of cells, seen from Cartesian coordinates, to be the same at the interface. That is, $h_r = r_* h_\theta = r_* h_s$ since $|e^{s_* \pm h_s} - e^{s_*}| \approx r_* h_s$. These meshes are shown in \Cref{fig:mesh}. We can see that the regular polar mesh includes rectangles of high aspect ratio near the origin, while the exponentially stretched one consists solely of square-shaped cells. In this sense, we regard exponential stretching as a natural treatment for polar coordinates.
    
    \begin{figure}[t]
        \centering
		\includegraphics[width = 180pt]{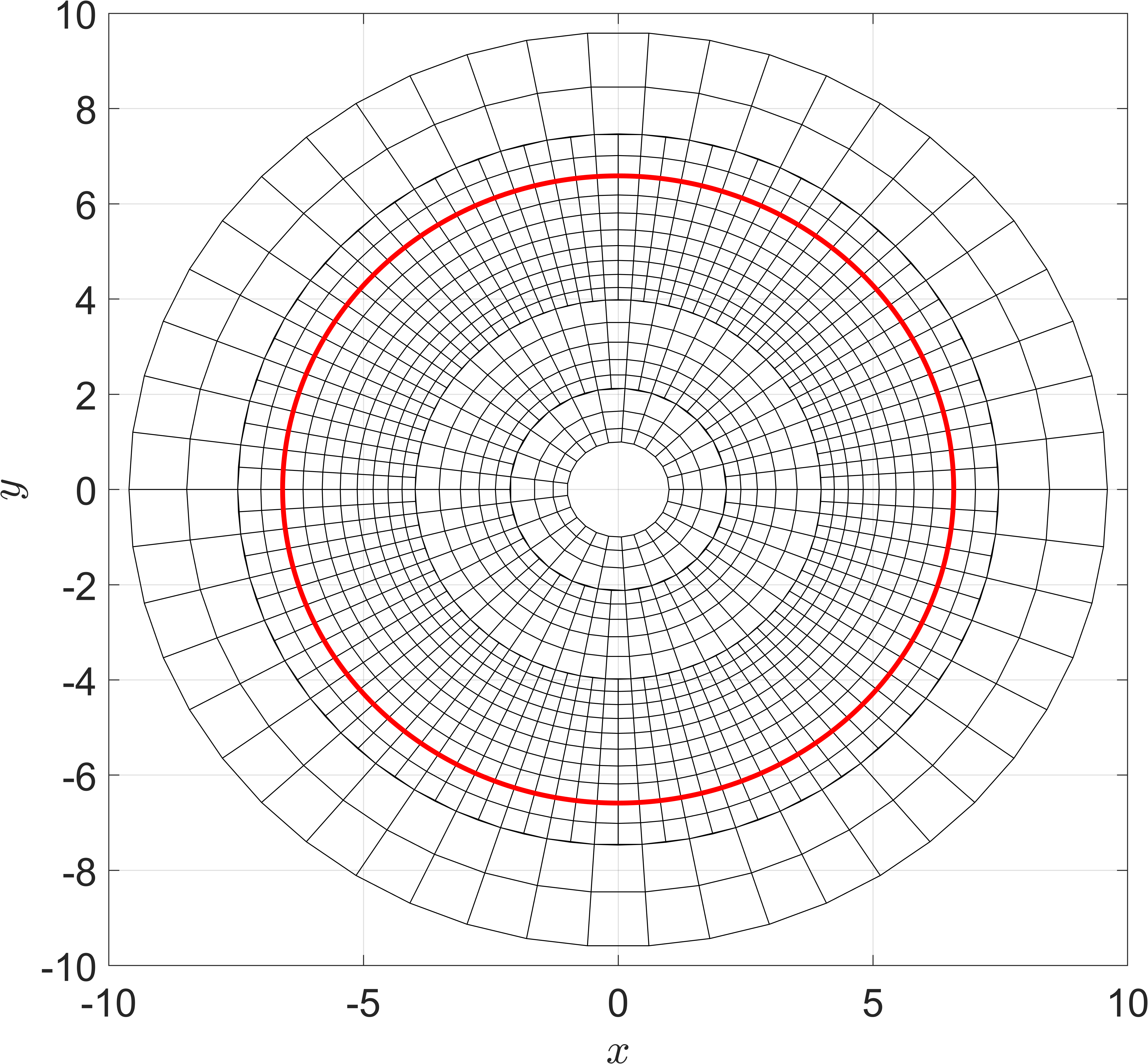}
		\includegraphics[width = 180pt]{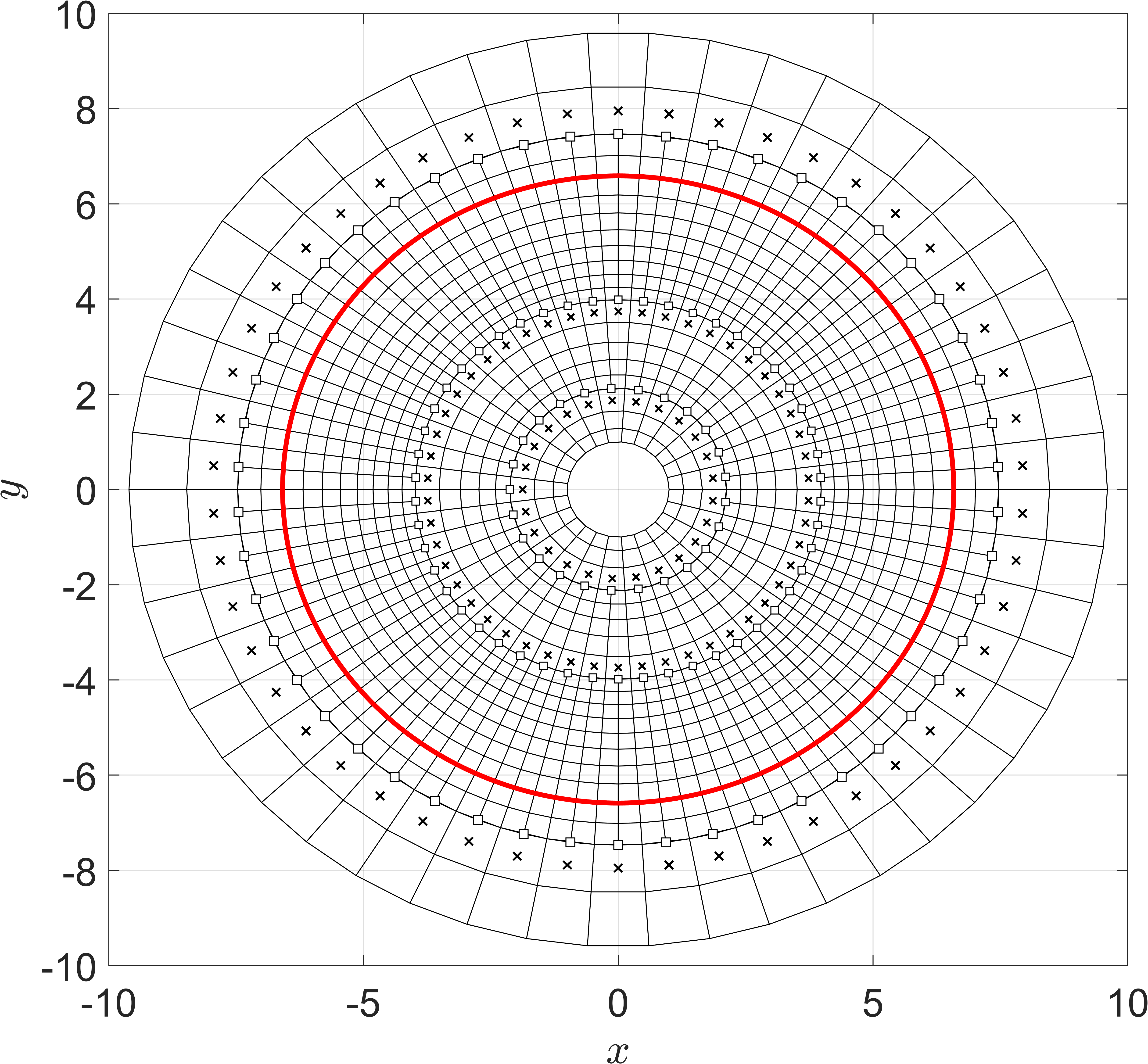}
        \caption{Result of mesh refinement under the exponentially stretched mesh. In the second figure, the dangling nodes and auxiliary nodes for the FDM are marked in empty squares and crosses, respectively.}
        \label{fig:refinement}
    \end{figure}
    
    On the other hand, the exponential stretching introduces cells at different scales, which can be resolved by a proper mesh refinement. An example is given in \Cref{fig:refinement}. We start with a coarse mesh near the scatterer and perform mesh refinement uniformly along $\theta$ direction as long as the $s$ coordinate increases by $\log 2$. This ensures the size of cells in the Cartesian coordinates differ by a factor of at most $2 + \bo(h)$. Combining mesh refinement with exponential stretching thus enables us to produce a quasi-uniform mesh in Cartesian coordinates. Furthermore, once we reach the PML layer, we can return to the coarse mesh as the solution and its derivatives decay sufficiently fast (\Cref{prop:sol_decay}). We can see that the mesh refinement effectively reduces the number of nodes under the same mesh size at the interface.
    
    Note that dangling nodes (empty squares in \Cref{fig:refinement}) appear when coarse cells transition to fine cells. Unfortunately, we are unable to find a reasonably small stencil that is centered at a dangling node and achieves sixth order consistency. Here, smallness means the stencil points are located within a short distance from the stencil center. To this end, we introduce auxiliary nodes (crosses in \Cref{fig:refinement}) located at the center of each coarse cell that is adjacent to a finer cell. These auxiliary nodes are then included in the stencil centered at dangling nodes. We will continue the discussion on stencils in the next section.
    
	\subsection{Modification to the FDM}
	\label{sec:modification}
	
	Since \cref{eq:PDE_PML_s} also falls within the general form of \eqref{eq:PDE_PML}, the methods discussed in \Cref{sec:FDM} are still applicable. Denote the set of nodes in the refined mesh to be $\Omega_h$. There are altogether three types of nodes in the refined mesh:
	\begin{itemize}
		\item[(a)] Interior nodes, where a compact stencil can be fit into $\Omega_h$;
		\item[(b)] Dangling nodes that appear at the border of the refinement region;
		\item[(c)] Auxiliary nodes at the center of coarse cells adjacent to a finer cell.
	\end{itemize}
	Besides, we need to reset the mesh size $h_s$ for each stencil locally. For an interior node $\spt$, $h_s$ is defined to be the smallest $h' > 0$ such that $\spt + h' \{-1, 0, 1\}^2 \in \Omega_h$. The other two types lies in a transition zone and $h_s$ is defined to be the size of the finer adjacent mesh.
	
\begin{figure}[th]
        \centering
		\includegraphics[height = 240pt]{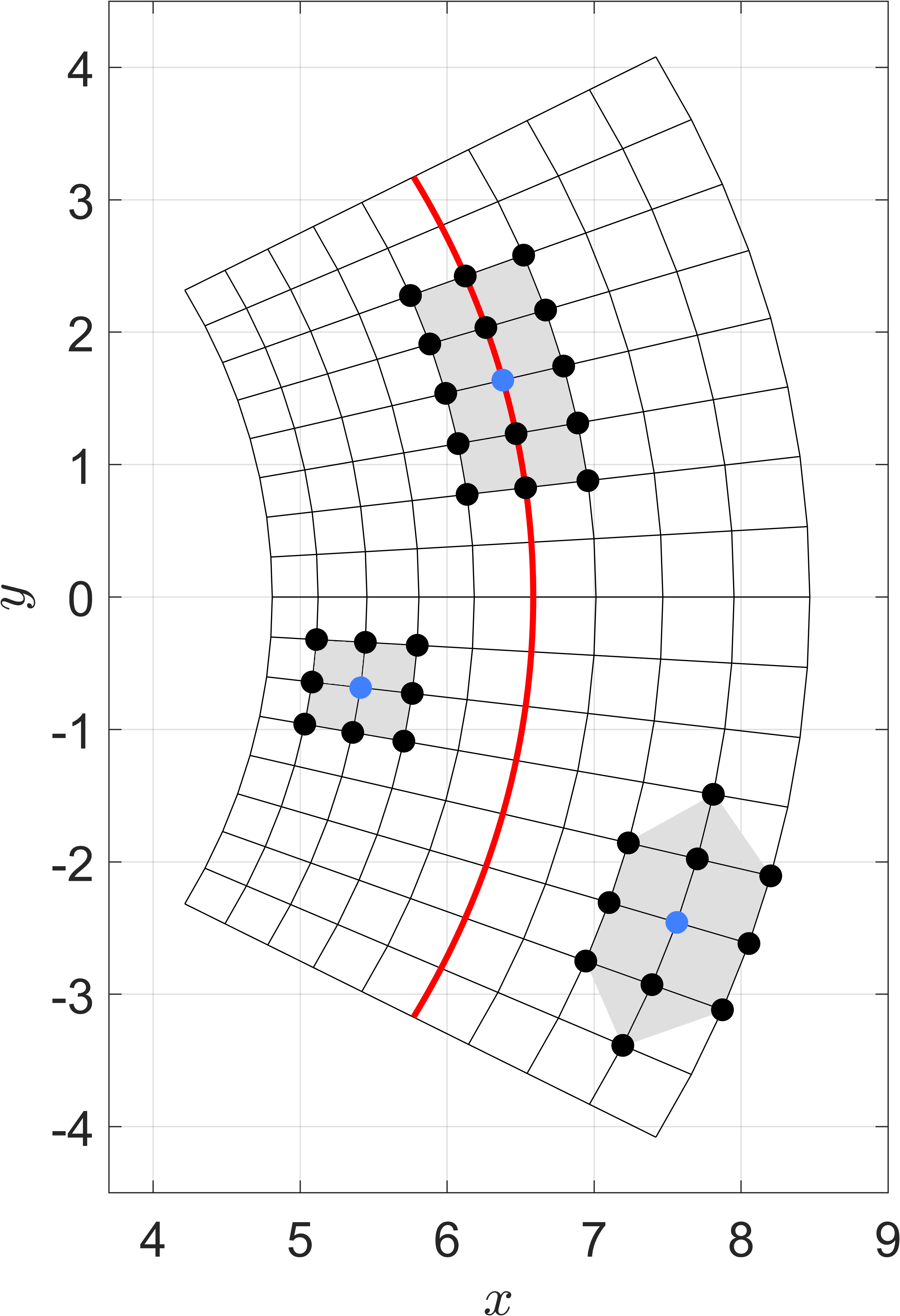}
		\includegraphics[height = 240pt]{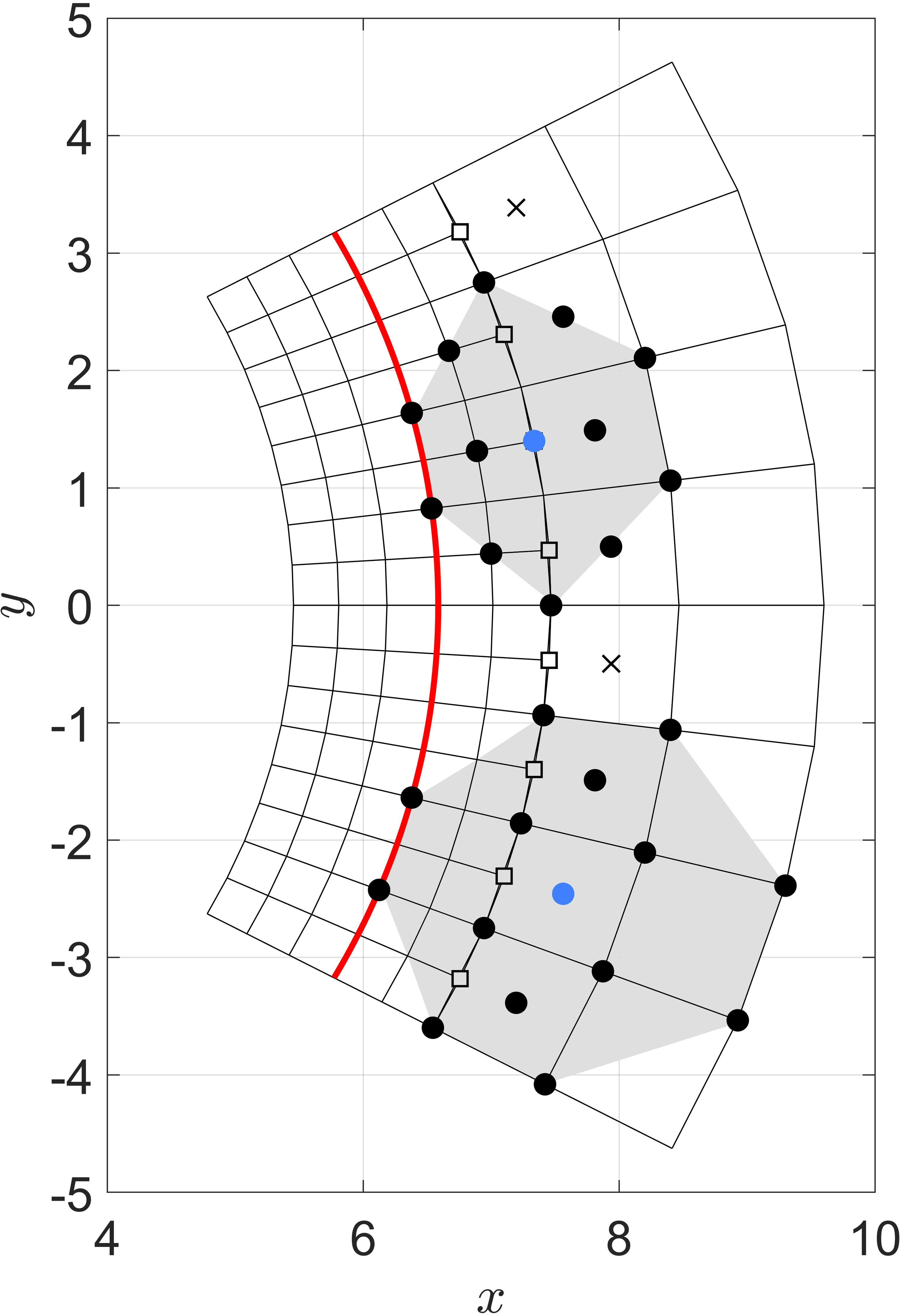}
        \caption{Different types of stencils with stencil center marked in blue and other stencil points marked in black. Stencils centered at an auxiliary node in $\Omega$ are not shown. The red curve, empty squares and crosses are the same as in \Cref{fig:refinement}.}
        \label{fig:stencil_s}
    \end{figure}

	As sixth order compact FDMs are available in the regular region $\Omega$, we suitably design the stencils elsewhere such that the overall consistency order is 6. \Cref{sec:mesh_descr} describes the default mesh refinement that is uniform in $\theta$ direction, yet we aim to perform mesh refinement freely so that our FDM can adapt to locally supported sources and non-circular scatterers. As there is no source and the solution decays uniformly in the PML region (see \Cref{prop:sol_decay}), we still assume that the uniform-in-$\theta$ refinement is performed in $\Gamma \cup \Omega^{\mathsf{PML}}$ for simplicity. Now we list all different types of stencils in \Cref{table:stencil_s} (also see \Cref{fig:stencil_s}). Here, the reference stencil $\SS$ is defined such that $\spt + h_s \SS$ is the stencil centered at $\spt$.
    
    \begin{table}[h]
        \small
        \begin{center}
        \begin{NiceTabular}{|p{48pt}|c||c|}[cell-space-limits=3pt]
            \hline
            \Block{1-2}{Location of stencil center $\spt$} & & Reference stencil $\SS$ \\ \hline \hline
            \Block{3-1}{In $\Omega$} & Interior & $\{-1, 0, 1\} \times \{-1, 0, 1\}$ \\ \cline{2-3}
            & \Block{1-1}{Dangling\\(facing $s$ direction)} & \Block{1-1}{$\big\{ (-2, -1), (-2, 1), (-1, -2), (-1, 0), (-1, 2), (0, -3)$,\\$(0, 0), (0, 3), (1, -2), (1, 0), (1, 2), (2, -1), (2, 1) \big\}$} \\ \cline{2-3}
            & \Block{1-1}{Dangling\\(facing $\theta$ direction)} & Exchange $s$ and $\theta$ components in the previous case \\ \cline{2-3}
            & Auxiliary & \Block{1-1}{$\big\{ (-3, -1), (-3, 1), (-1, -3), (-1, -1), (-1, 1), (-1, 3)$,\\$(0, 0), (1, -3), (1, -1), (1, 1), (1, 3), (3, -1), (3, 1) \big\}$} \\ \hline
            \Block{1-1}{On $\Gamma$} & Interior & $\{-1, 0, 1\} \times \{-2, -1, 0, 1, 2\}$ \\ \hline
            \Block{3-1}{In $\Omega^{\mathsf{PML}}$} & Interior & $\{-1, 0, 1\} \times \{-1, 0, 1\} \cup \{ (0, \pm 2) \}$ \\ \cline{2-3}
            & \Block{1-1}{Dangling\\(facing $s$ direction)} & Same as dangling nodes in $\Omega$ \\ \cline{2-3}
            & Auxiliary & \Block{1-1}{$\SS \cup \{ (0, \pm 2) \}$, where $\SS$ corresponds to\\the stencil at an auxiliary node in $\Omega$} \\ \hline
        \end{NiceTabular}
        \end{center}
        \caption{All types of reference stencils $\SS$ in $(s, \theta)$ coordinates.}
        \label{table:stencil_s}
    \end{table}
    
	A few remarks are in order. First, except for the center, the stencils at a dangling node only consist of interior and auxiliary nodes, and the stencils at an auxiliary node only consist of interior nodes from the coarser mesh. This shows that the FDM is well-defined. Moreover, for stencils centered on dangling and auxiliary nodes, the maximum distance of stencil points from the center only slightly increases from $\sqrt{8} h_s$ to $\sqrt{10} h_s$, which maintains the local truncation error almost at the same level. There are some minor restrictions when we implement mesh refinement. For example, there is no consecutive refinements within two layers of nodes. We will not elaborate on such details.
	
	Now we wish to apply \Cref{thm:pollution} as in the case of regular polar coordinates. The zeroth order stencil coefficients in the regular region $\Omega$ are listed below:
	
	\vspace{0.5em}%
	\noindent Interior: $c_{(\pm1, \pm1), 0} = -1$, $c_{(\pm1, 0), 0} = c_{(0, \pm1), 0} = -4$, $c_{(0, 0), 0} = 20$.
	
	\noindent Dangling (facing $s$ direction): $c_{(\pm2, \pm1), 0} = -33$, $c_{(\pm1, \pm2), 0} = -75$, $c_{(\pm1, 0), 0} = -378$, $c_{(0, \pm3), 0} = -14$, $c_{(0, 0), 0} = 1216$.
		
	\noindent Dangling (facing $\theta$ direction): $c_{(\pm1, \pm2), 0} = -33$, $c_{(\pm2, \pm1), 0} = -75$, $c_{(0, \pm1), 0} = -378$, $c_{(\pm3, 0), 0} = -14$, $c_{(0, 0), 0} = 1216$.
	
	\noindent Auxiliary: $c_{(\pm3, \pm1), 0} = c_{(\pm1, \pm3), 0} = -1$, $c_{(\pm1, \pm1), 0} = -14$, $c_{(0, 0), 0} = 64$.
	
	\vspace{0.5em}%
	\noindent Through symbolic computation, we have also verified that these stencil coefficients are unique up to constant multiples. Hence, we can perform Algorithm 1 with $\delta_h = 0$ within the regular region and $\delta_h \sim h^{16}$ elsewhere to obtain a sixth order FDM with reduced pollution effect.
	
    \subsection{Handling non-circular scatterers}
    \label{sec:non_circular}
    
    Consider a node $\spt \in \Omega_h$ near the scatterer, which we will refer to as a boundary node. A stencil centered at (not at) a boundary node is referred to as a boundary (non-boundary) stencil. Since $\spt \in \Omega$, we only need to treat the equation $\Delta v + \kappa^2 e^{2s} v = e^{2s} f$ and the boundary condition $v = g$ near $\spt$. We set $\tilde{f} = e^{2s} f$. 
    
    A good boundary stencil is crucial to the stability and eventual convergence of FDMs. In view of the isotropy of Laplace operator as opposed to \eqref{eq:PDE_PML_r}, we adapt the techniques in \cite{han2025convergent}, which proposes a convergent sixth order FDM in curved domains. Previously discussed mesh refinement further enhances our FDM as it can fit into various scatterer geometries. This section is divided into two parts: a summary of the generic method in \cite{han2025convergent}, and our pollution minimization technique for the boundary stencils. The generic method is introduced just for the sake of completeness, and consists of procedures unnecessary for the pollution minimization technique. From the view of implementation, one can completely ignore the contents from \cref{eq:v_taylor_bdr0} (excluding itself) to the end of its next paragraph.
    
	The argument in \cite{han2025convergent} is in parallel to \cite{feng2025symmetric}, which we have applied to treat non-boundary stencils. However, \cite{han2025convergent} involves complex partial derivatives $\partial_{\C}^{\bm{k}} v := 2^{-|\bm{k}|} (\partial_s - \ii \partial_\theta)^{k_1} (\partial_s + \ii \partial_\theta)^{k_2} v$ instead of regular partial derivatives to make full use of the boundary condition. The authors extensively used the fact that the PDE and its solution are real-valued in order to simplify complex-valued expressions. When applying the formulas in \cite{han2025convergent} into our situation of $v \in \C$, one needs to regard $\overline{\partial_{\C}^{(\ell_1, \ell_2)}} v$ as $\partial_{\C}^{(\ell_2, \ell_1)} v$ instead of $\partial_{\C}^{(\ell_2, \ell_1)} \overline{v}$. The same is true for real and imaginary parts of an expression. To have a clearer comparison between our setting and the setting in \cite{han2025convergent}, we temporarily treat $v$ as if it is real-valued in this section. Again, we only mention certain essential steps and the differences compared to \cite{han2025convergent}.
    
    Analogous to \eqref{eq:v_taylor_int}, we can expand the solution $v$ in terms of complex partial derivatives:
	\begin{equation}
		\label{eq:v_taylor_int_complex}
        v(\bpt + ph)
        = \sum_{\bm{\ell} \in \LNC_{M - 1}} \sum_{k = |\bm{\ell}|}^{M - 1}
        A^k_{\bm{\ell}, \C} (p) h_s^k \cdot \partial_{\C}^{\bm{\ell}} v(\bpt) + F_{\C}(p) + \bo (h_s^M).
	\end{equation}
	Here $\bpt \in \partial D$ is a point close to, but not the same as the stencil center $\spt$, $\LNC_{M + 1} := \{ \bm{k} = (k_1, k_2) \in \N_0^2: |\bm{k}| \leq M + 1, k_1 k_2 = 0 \}$, and $A^k_{\bm{\ell}, \C} (p)$ and $F_{\C}(p)$ are explicitly known quantities involving the derivatives of $\kappa^2 e^{2s}$ and $\tilde{f} = e^{2s} f$ at $\bpt$. Note that we have set $M$ lower by 2 than before since \cite{han2025convergent} indicates that this is enough to achieve the same accuracy order for boundary stencils. The differences are: (1) $\tilde{a}^{\bm{k}}_{\bm{\ell}} = -\frac{1}{4} \binom{k_1 - 1}{\ell_1} \binom{k_2 - 1}{\ell_2} \kappa^2 e^{2s}$ instead of \cite[equation (2.9)]{han2025convergent}, and (2) the simplified version of \eqref{eq:v_taylor_int_complex} (i.e., \cite[equation (2.19)]{han2025convergent}) needs to be modified into
	\begin{equation}
		\label{eq:v_taylor_bdr0}
		v(\bpt + ph) = \sum_{m = 0}^{M - 1} \sum_{k = m}^{M - 1} 2 \re \left( (1 - \tfrac{1}{2} \td(m)) A^k_{(m, 0), \C}(p) \partial_{\C}^{(m, 0)} v(\bpt) \right) h_s^k + F_{\C}(p) + \bo(h_s^M),
	\end{equation}
	as $A^k_{(0, 0), \C} = \td(k)$ no longer holds. One can equivalently replace $\tilde{A}^{\bm{k}}_{\bm{\ell}}$ in \cite[equation (2.11)]{han2025convergent} by $(1 - \frac{1}{2} \td(\bm{\ell})) \tilde{A}^{\bm{k}}_{\bm{\ell}}$ and keep the remaining formulas unchanged to derive
	\begin{equation*}
		v(\bpt + ph) = \sum_{m = 0}^{M - 1} \sum_{k = m}^{M - 1} 2 \re \left( A^k_{(m, 0), \C}(p) \partial_{\C}^{(m, 0)} v(\bpt) \right) h_s^k + F_{\C}(p) + \bo(h_s^M).
	\end{equation*}
	
	By utilizing the boundary conditions, we can obtain a better expansion of $v$ with fewer unknown complex partial derivatives: 
	\begin{equation}
		\label{eq:v_taylor_bdr}
        v(\bpt + ph)
        = \sum_{m = 1}^{M - 1} \sum_{k = m}^{M - 1}
        A^k_{m, \C} (p) h_s^k \cdot \im (e^{\ii m \theta} \partial_{\C}^{(m, 0)} v(\bpt)) + G_{\C}(p) + \bo (h_s^M).
	\end{equation}
	Here $\theta$ is the tangent angle at $\bpt$ and the calculation of various quantities $A^k_{m, \C}(p)$ and $G_{\C}(p)$ is given in \cite{han2025convergent}. We summarize the differences below. Suppose the quantity $\tilde{A}^{\bm{k}}_{\bm{\ell}}$ has already been replaced as above, then all the proof and calculation in pages 10-11 of \cite{han2025convergent} remain the same, except that: (1) the summation indices $m$ and $n$ begins from $0$; (2) the range of $m$ and $n$ become lower by $1$ in \cite[equations (4.7)-(4.9)]{han2025convergent}; (3) all $u(\bpt)$ and $g(\bpt)$ are absorbed into the summation and hence disappear from the formulas. After we obtained the expansion \eqref{eq:v_taylor_bdr}, we proceed the same way as in \cite{han2025convergent} to construct an FDM. In particular, the determination of boundary stencils are described in Section 4.2 and the stencil coefficients are obtained using Lemma 4.2 in that article.
	
	Now we consider the pollution minimization technique for boundary stencils. As it is difficult to find a test function $w$ that satisfies both $\Delta w + \kappa^2 w = 0$ and $w = 0$ on a boundary segment, the calculation of $G_{\C}(p)$ for the test functions $w(\cdot; \theta_0)$, $\theta_0 \in \T$ must be included. The very complicated expression of $G_{\C} (p)$ makes the pollution minimization technique both computationally ineffective and unstable. Hence, in what follows, we essentially treat a boundary stencil as a special interior stencil instead. Note that we still replace $M$ by $M - 2$ in this situation.
	
	Using complex partial derivatives, the stencil coefficients $C_p (\spt) = \sum_{j = 0}^{M - 1} c_{p, j} h^j$ for an interior stencil satisfies	
	\begin{equation}
		\label{eq:c_p_j_complex}
		\sum_{p \in \SS} A^{|\bm{\ell}|}_{\bm{\ell}, \C}(p) c_{p, j}
		= -\sum_{k = 0}^{j - 1} \sum_{p \in \SS} A^{|\bm{\ell}| + j - k}_{\bm{\ell}, \C}(p) c_{p, k}, \quad \forall \, \bm{\ell} \in \LNC_{M - 1 - j}
	\end{equation}
	according to \cite[Lemma 3.1]{han2025convergent}. Here $\SS$ is defined such that $\bpt + h \SS$ is the stencil centered at $\spt$, and $\bpt$ is the associated base point of $\spt$ in \eqref{eq:v_taylor_bdr0}. For any $0 \leq j \leq M - 1$, \eqref{eq:c_p_j_complex} is a linear system with at most $2M - 1$ equations. Setting $M = 6$, this linear system suggests that we need a stencil consisting of at least $12$ points. This should probably guarantee the existence of a nontrivial solution for the stencil coefficients except for the zero-measure event of the wide matrix $( A^{|\bm{\ell}|}_{\bm{\ell}, \C}(p))_{\bm{\ell} \in \LN_{M - 1 - j}, p \in \SS}$ not having full row rank. Bearing the stability of the FDM in mind, we start with the suggested stencil in \cite[Section 4.2]{han2025convergent} and enlarge it with points on the boundary. Since the stencils in \cite{han2025convergent} consist of $7$ or $8$ points, we choose $5$ extra points from $\partial \Omega$ that are of distance $\sim h$ with each other. For example, suppose $\partial \Omega$ locally consists of a curve $\{ (\gamma_x(t), \gamma_y(t)): t \in -(\delta, \delta) \}$ with $(\gamma_x(0), \gamma_y(0)) = \bpt$ and denote $\sigma = \big( \gamma_x'(0)^2 + \gamma_y'(0)^2 \big)^{1/2}$. The additional $5$ points can be taken as $\{ (\gamma_x(k/\sigma), \gamma_y(k/\sigma)): k = -2, -1, 0, 1, 2 \}$. After we have fixed the stencil, we find the stencil coefficients using the pollution minimization process exactly the same as in \Cref{sec:pollution} (instead of solving \eqref{eq:c_p_j_complex}). Finally, we set up the FDM $\mathcal{L}_h v_h (\spt) = \sum_{p \in \SS} C_p (\spt) F_{\C}(p)$, which is expected to have sixth order consistency in terms of boundary stencils:
	\begin{equation*}
		\mathcal{L}_h v (\spt) = \sum_{p \in \SS} C_p (\spt) F_{\C}(p) + \bo(h^6)
		\quad \text{for all} \ v \ \text{satisfying} \ \Delta v + \kappa^2 e^{2s} v = \tilde{f}.
	\end{equation*}
	The value of $v_h$ at stencil points on the true boundary can be directly obtained from the boundary condition and moved to the right hand side of the FDM linear system.
    
	\section{Numerical results}
	\label{sec:numerical}
	
	\subsection{Choice of parameters and general comments}
	\label{sec:choice_param}
	
	First, we briefly talk about the position of the PML layer and the generation of meshes before refinement. We predetermine two values $\tilde{r}_*$, $\tilde{r}_M$ and later adjust the PML layer to be close to $\{\tilde{r}_* < r < \tilde{r}_M\}$. As we have discussed in \Cref{sec:mesh_descr}, we need to let $h_r = r_* h_\theta = r_* h_s$ so that the mesh size at the interface is the same from Cartesian coordinates. In practice, we set $h := h_r = \tilde{r}_* h_\theta = \tilde{r}_* h_s$ and refer to this value as the mesh size at the interface. We will take $h = \frac{2\pi \tilde{r}_*}{N}$ with $N \in \N$ so that the nodes are evenly distributed in $\theta$ direction. We then adjust $r_*$, $r_M$, $s_*$, $s_M$ close to $\tilde{r}_*$, $\tilde{r}_M$, $\log \tilde{r}_*$ and $\log \tilde{r}_M$ respectively to align $\partial D$ and $\partial \Omega^{\mathrm{PML}}$ with grid lines. Finally, unless otherwise stated, we adopt the way in \Cref{sec:mesh_descr} to perform mesh refinement.
	
	We next discuss some topics related to the complex coordinate transform $\rho$. In the examples we will take $\rho(r) = r + \ii \alpha_1 (r - r_*)^2$ for the regular polar coordinates and $\rho(s) = s_* + \alpha_2 (s - s_*)$ for the exponentially stretched polar coordinates with $\alpha_1 \in \R$ and $\alpha_2 \in \C$. Let $d := r_M - r_* = e^{s_M} - e^{s_*}$ be the thickness of the PML layer. Both \Cref{ex:ex2} and \cite[Section 5.1]{yang2021truly} indicate that $d = \bo(\kappa^{-1})$ can ensure a good performance of the numerical methods. Since the mesh size should be no greater than $\bo(\kappa^{-1})$ as well to resolve a wave with wavenumber $\kappa$, we see that the PML layer consists of a constant number of nodes in the radial direction. On the contrary, if we impose the factor $\kappa^{-1}$ on the function $\rho$, then the number of nodes in the PML layer will significantly increase, resulting in a loss of efficiency. According to the decay estimate \eqref{eq:sol_decay}, we pick $\alpha_1$, $\alpha_2$ such that $\exp(-\kappa \im z (1 - r_*^2 / |z|^2)^{1/2}) \lesssim 10^{-16}$ at $z = \rho(r_M)$ or $z = e^{\rho(s_M)}$. More explicitly, we determine $\alpha_1$ and $\alpha_2$ as follows:
	\begin{equation*}
		\alpha_1 = 40 \kappa^{-1} d^{-2} (1 - r_*^2 / r_M^2)^{-1/2}, \quad
		\alpha_2 = \left( \tfrac{1}{2} t^{-1} \log(1 + (40 \kappa^{-1} e^{-s_*})^2) + \ii \arcsin(t) \right) (s_M - s_*)^{-1}
	\end{equation*}
	with $t \in (0, 1]$ as a tunable parameter.
	
	The regularization parameter $\delta_h$ in the pollution minimization process is determined as follows. We take $\delta_h = 0$ for non-boundary stencils in the regular region. For all other stencils, we take $\delta_h = \max \{ 0.01 \widetilde{\mathcal{I}}_h (\vec{a}_*), 10^{-14} \}$, where $\widetilde{\mathcal{I}}_h$ is defined in \eqref{eq:I_h} and $\vec{a}_*$ is computed by minimizing $\widetilde{\mathcal{I}}_h$ instead of $\mathcal{I}_h$ in \eqref{eq:hat_C_p}. This ensures $\delta_h$ decreases in the speed of $h^{2M + 4}$ (or $h^{2M}$ on the boundary) when $h$ is not too small. For the case of small $h$, we follow \Cref{app:float_error} to set a threshold of $10^{-14}$ to $\delta_h$. This stabilizes the local truncation error is at the magnitude of $10^{-7}$. 
	
	Let $v_h$ be the numerical solution to \eqref{eq:PDE_PML_r} or \eqref{eq:PDE_PML_s}. If the exact solution $u$ to \eqref{eq:PDE} is available, we take the error $e_h$ as
	\begin{equation*}
		e_h = \| v_h - u \|_{\infty} = \max_{\spt \in \Omega_h} |v_h(\spt) - u(\spt)|.
	\end{equation*}
	Otherwise, we generate a reference solution $v^{\mathrm{ref}}$ using a fine mesh $\Omega_h^{\mathrm{ref}}$, and determine the error by
	\begin{equation*}
		e_h = \| v_h - v^{\mathrm{ref}} \|_{\infty} = \max_{\spt \in \Omega_h \cap \Omega_h^{\mathrm{ref}}} |v_h(\spt) - v^{\mathrm{ref}}(\spt)|.
	\end{equation*}
	Note that if $h / h^{\mathrm{ref}} \in \N$ with $h^{\mathrm{ref}}$ being the mesh size of the reference solution, then $\Omega_h$ will almost be a subset of $\Omega_h^{\mathrm{ref}}$.
	
	\subsection{Examples}
	\label{sec:examples}
	
	\begin{figure}[b]
    	\centering
		\includegraphics[width = 0.35 \linewidth]{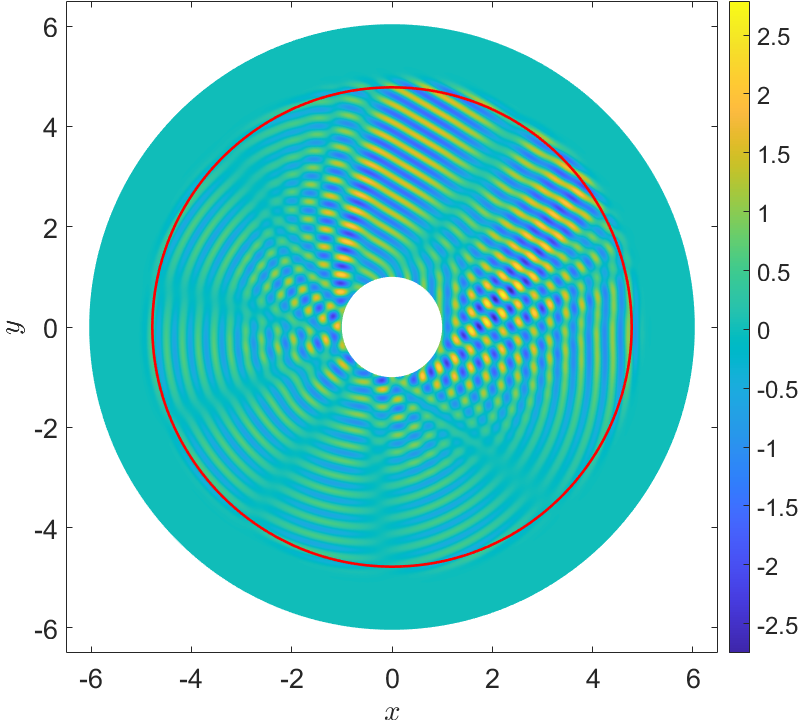}
		\hspace{0.02 \linewidth}
		\includegraphics[width = 0.35 \linewidth]{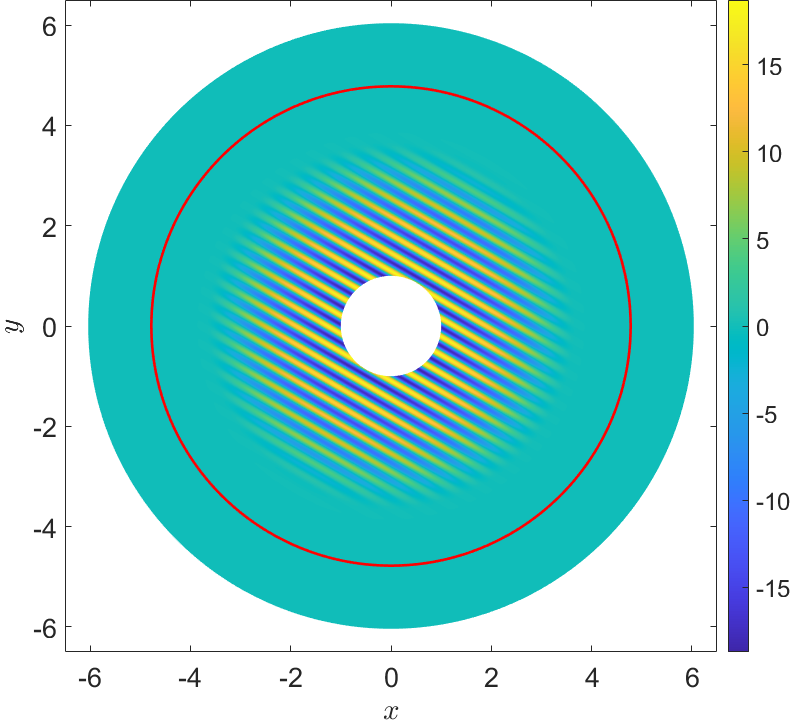}
		\vspace{-0.5em}
       	\caption{Real part of reference solution $v^{\mathrm{ref}}$ (left) and source term $f$ (right) in \Cref{ex:ex1}. Red circle indicates the interface.}
	    \label{fig:ex1:vf}
	\end{figure}
	
	\begin{figure}[b]
	\centering
    \begin{minipage}[t]{0.7 \linewidth}
    	\vspace{0pt}
    	\centering
        \begin{NiceTabular}{|c||c|c||c|c||c|c|}[cell-space-limits=2pt]
            \hline
            & \Block{1-2}{Regular polar\\coordinates} & & \Block{1-2}{Exponential\\stretching} & & \Block{1-2}{With mesh\\refinement} \\ \hline
            $\kappa h$ & $e_h$ & order & $e_h$      & order & $e_h$      & order \\ \hline \hline
            2.099 & 4.432E$-$1 &       & 1.114E$-$2 &       & 3.070E$-$1 &       \\ \hline
            1.574 & 1.081E$-$1 & 4.91  & 1.877E$-$3 & 6.19  & 7.379E$-$3 & 12.96 \\ \hline
            1.049 & 1.553E$-$2 & 4.78  & 1.509E$-$4 & 6.22  & 4.617E$-$4 & 6.84  \\ \hline
            0.787 & 4.977E$-$3 & 3.96  & 2.456E$-$5 & 6.31  & 7.767E$-$5 & 6.20  \\ \hline
            0.525 & 8.903E$-$4 & 4.24  & 2.909E$-$6 & 5.26  & 6.145E$-$6 & 6.26  \\ \hline
            0.262 & \Block{1-6}{Reference solution} & & & & & \\ \hline
        \end{NiceTabular}
    \end{minipage}
    \begin{minipage}[t]{0.23 \linewidth}
    	\vspace{-0.9em}
        \centering
        \includegraphics[width = \linewidth]{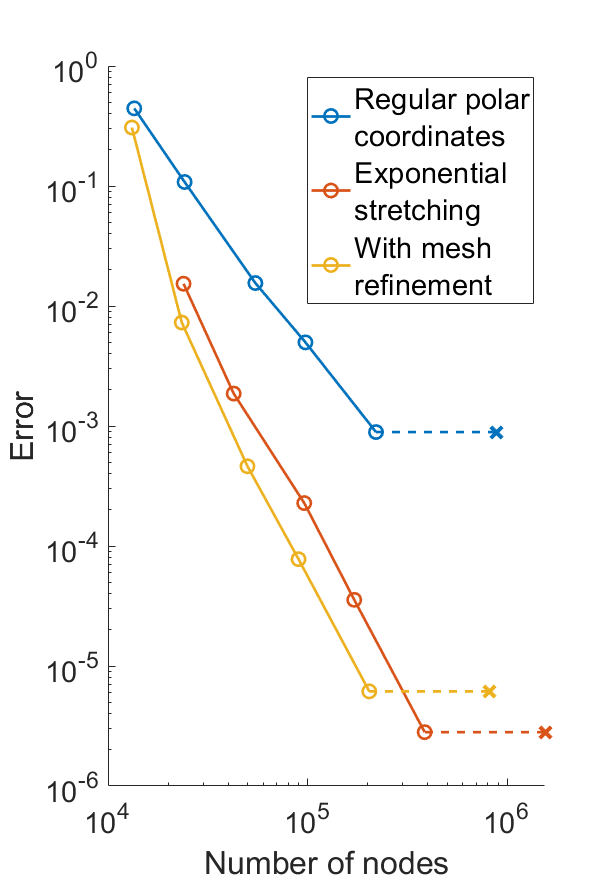}
    \end{minipage}
    
	\captionof{table}{Error and convergence order using three different meshes in \Cref{ex:ex1}. The mesh size $h$ at the interface is the same for all data on the same row. The order corresponding to mesh size $h$ is calculated by $\log (e_{h'} / e_h) / \log (h' / h)$, where $h'$ is the previous mesh size.}
	\label{table:ex1:order}
	\vspace{-0.5em} \captionof{figure}{Comparison of numerical error versus number of nodes with different meshes in \Cref{ex:ex1}. Data points in the same ordinal position (first, second, etc.) across three lines use the same mesh size $h$ at the interface. The rightmost data points correspond to the reference solutions, and only the number of nodes matters.}
	\label{fig:ex1:nnodes}
    \end{figure}
	
	\begin{example}
		\label{ex:ex1} \normalfont
		In this example, we compare the performance and verify the accuracy order using different meshes: regular polar mesh, and exponentially stretched mesh with and without mesh refinement. The parameters and functions in \eqref{eq:PDE} are given by $\kappa = 20$, $D = \{ (r, \theta): r < 1 \}$ and 
		\begin{align*}
			& f(r, \theta) = \kappa \exp \big( \ii \kappa r \cos(\tfrac{\pi}{3} - \theta) \big) \cdot \boldone_{\{ r < 4 \}} \exp \big( \tfrac{r^2}{r^2 - 16} \big), \\
			& g(\theta) = u|_{\partial D} (\theta) = \exp(\ii \kappa \cos \theta). 
		\end{align*}
		Note that the source and boundary data are misaligned, and the exact solution $u$ is unknown. The PML layer is adjusted from $\tilde{r}_* = e^{0.5 \pi}$ and $\tilde{r}_M = 1.25 \tilde{r}_*$, which results in $\kappa d \approx 25.2$ after the mesh is set. We take $h^{\mathrm{ref}} = \frac{2\pi \tilde{r}_*}{2304}$ and calculate numerical solutions using $h = N h^{\mathrm{ref}}$, $N = 1, 2, 3, 4, 6, 8$. 
		
		We present the reference solution $v^{\mathrm{ref}}$ and the source term $f$ in \Cref{fig:ex1:vf}. \Cref{table:ex1:order} lists the error and convergence order using the same mesh size. This shows that our FDMs achieve the anticipated accuracy order. We further compare these three meshes featuring the number of nodes $\# \Omega_h$ and the error $e_h$ in \Cref{fig:ex1:nnodes}. We can see that mesh refinement successfully reduces the number of nodes to the level of regular polar coordinates while maintaining the improved accuracy from exponential stretching. Finally, \Cref{table:ex1:pollution} verifies the superior performance of our pollution minimization technique which is able to reduce numerical error by up to 99.78\%. The FDM without pollution minimization is constructed based on \Cref{sec:FDM_generic}.
	\end{example}

    \begin{table}[t]
    	\centering
        \begin{NiceTabular}{|c||c|c||c|c||c|c|}[cell-space-limits=2pt]
            \hline
            & \Block{1-2}{Regular polar\\coordinates} & & \Block{1-2}{Exponential\\stretching} & & \Block{1-2}{With mesh\\refinement} \\ \hline
            $\kappa h$ & $e_h'$ & $R$     & $e_h'$     & $R$     & $e_h'$     & $R$     \\ \hline \hline
            2.099  & 3.956E$+$0 & 88.80\% & 4.744E$+$0 & 99.77\% & 9.084E$+$0 & 96.62\% \\ \hline
            1.574  & 1.954E$+$0 & 94.47\% & 6.727E$-$1 & 99.72\% & 1.762E$+$0 & 99.58\% \\ \hline
            1.049  & 4.679E$-$1 & 96.68\% & 6.109E$-$2 & 99.75\% & 1.020E$-$1 & 99.55\% \\ \hline
            0.787  & 1.385E$-$1 & 96.41\% & 1.123E$-$2 & 99.78\% & 1.651E$-$2 & 99.53\% \\ \hline
            0.525  & 2.102E$-$2 & 95.76\% & 8.931E$-$4 & 99.67\% & 1.356E$-$3 & 99.55\% \\ \hline
            0.262  & \Block{1-6}{Reference solution} & & & & & \\ \hline
        \end{NiceTabular}
        \caption{Numerical error $e_h'$ from FDMs without pollution minimization process and the percentage $R$ of error reduced by the pollution minimization process. $R$ is calculated by $R := 1 - e_h / e_h'$, where $e_h$ is the error in \Cref{table:ex1:order}.}
       	\label{table:ex1:pollution}
    \end{table}
    
	\begin{figure}[b]
		\vspace{-1em}
	    \centering
		\includegraphics[width = 0.4 \linewidth]{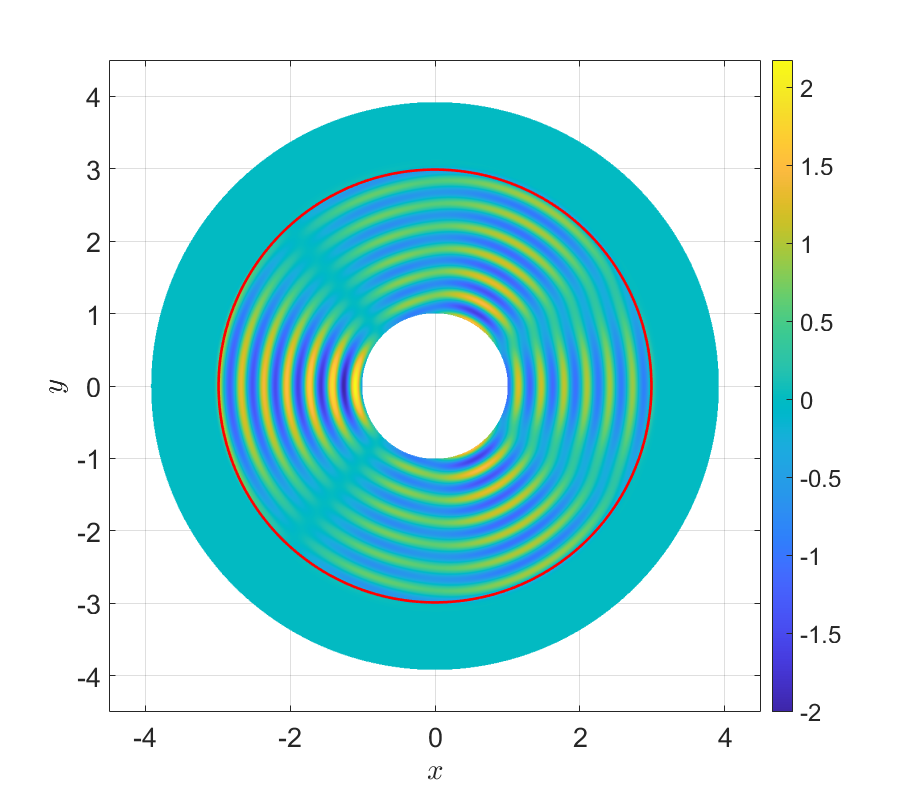}
		\includegraphics[width = 0.4 \linewidth]{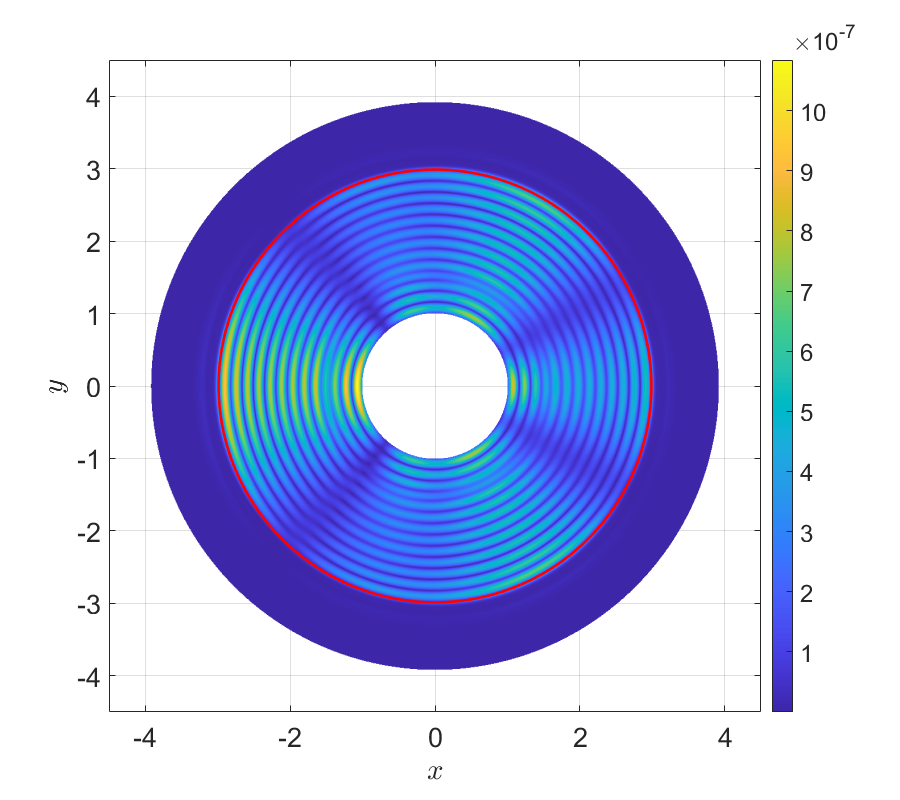}
        \caption{Real part of exact solution $u$ to \Cref{eq:PDE} (left) and numerical error $|u - u_h|$ (right) in \Cref{ex:ex2} with $\kappa = 20$ and $\kappa \tilde{d} = 20$. Red circle indicates the interface.}
        \label{fig:ex2:sol}
	\end{figure}
    
	\begin{example}
		\label{ex:ex2} \normalfont
		
		This example aims to show that our FDM is robust under a wide range of wavenumber $\kappa$ and PML thickness $d$. The parameters and functions in \eqref{eq:PDE} are given by $\kappa = 5$, $20$, $50$, $100$, $D = \{ (r, \theta): r < 1 \}$, $f = 0$ and 
		\begin{equation*}
			g(\theta) = u|_{\partial D} (\theta) = \sum_{j \in \Z} j^2 e^{-|j|} e^{\ii j (j + \theta)}.
		\end{equation*}
		
		\vspace{-0.5em}\noindent The solution $u$ to \eqref{eq:PDE} is known and given by
		\begin{equation*}
			u(\rho(r), \theta) = \sum_{j \in \Z} j^2 e^{-|j|} \frac{J_j (\kappa)}{H^{(1)}_j (\kappa)} H^{(1)}_j (\kappa \rho(r)) e^{\ii j (j + \theta)}.
		\end{equation*}
		The PML layer is adjusted from $\tilde{r}_* = 3$ and $\tilde{r}_M = \tilde{r}_* + \tilde{d}$ with $\tilde{d}$ satisfying $\kappa \tilde{d} = 3$, $5$, $10$, $20$, $30$, $50$. To better showcase the ability of the proposed FDM, we will use a moderately sized mesh with $\kappa h \approx 0.5$ for all $\kappa$ values. In view of its superior performance in \Cref{ex:ex1}, we adopt the exponentially stretched mesh with mesh refinement and apply the pollution minimization technique.
		
	\begin{figure}[t]
		\centering
		\begin{minipage}[t]{0.54 \linewidth}
	    	\vspace{0pt}
	        \begin{center}
	        \begin{NiceTabular}{|c||c|c|c|c|c|}[cell-space-limits=3pt]
	            \hline
	            $\kappa \tilde{d} \ \bs \ \kappa$ & 5 & 20 & 50 & 100  \\ \hline \hline
	            3  & 2.25E-3 & 8.54E-3 & 8.60E-4 & 4.44E-6 \\ \hline
	            5  & 9.62E-5 & 9.53E-4 & 4.97E-5 & 2.06E-6 \\ \hline
	            10 & 1.69E-5 & 4.06E-5 & 1.19E-6 & 1.40E-6 \\ \hline      
	            20 & 1.26E-6 & 1.09E-6 & 1.14E-6 & 8.04E-7 \\ \hline
	            30 & 4.63E-7 & 3.71E-7 & 1.12E-6 & 7.82E-7 \\ \hline
	            50 & 2.55E-7 & 3.85E-7 & 1.15E-6 & 8.00E-7 \\ \hline
	        \end{NiceTabular}
		    \end{center}
		\end{minipage}
		\begin{minipage}[t]{0.32 \linewidth}
	    	\vspace{-1em}
		    \centering
		    \includegraphics[width = \linewidth]{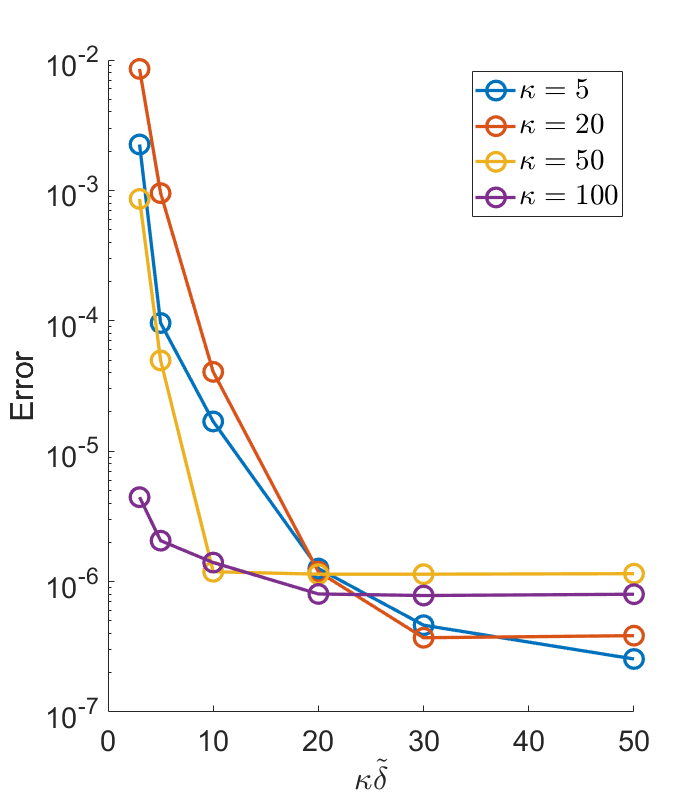}
		\end{minipage}
    
		\caption{Numerical error $e_h$ with different wavenumber $\kappa$ and (approximate) PML layer thickness $\tilde{d}$ in \Cref{ex:ex2}.}
		\label{fig:ex2:err}
		\vspace{-0.5em}
    \end{figure}
		
		\Cref{fig:ex2:sol} shows the exact solution and the numerical error with $\kappa = 20$, $\kappa \tilde{d} = 20$. In \Cref{fig:ex2:err} we present the numerical error over different values of $\kappa$ and $\tilde{d}$ mentioned above. We can see that the numerical error stabilizes at a very low level for all $\kappa$ values as long as $\kappa \tilde{d}$ is over a threshold of $20 \sim 30$. Even when the PML layer is as thin as $3\kappa^{-1}$ to $5\kappa^{-1}$, the FDM provides a decent error of magnitude $10^{-3}$. This hence confirms the robustness of our FDM.
	\end{example}
	
	\begin{figure}[b]
		\vspace{-1em}
	    \centering
		\includegraphics[width = 0.3 \linewidth]{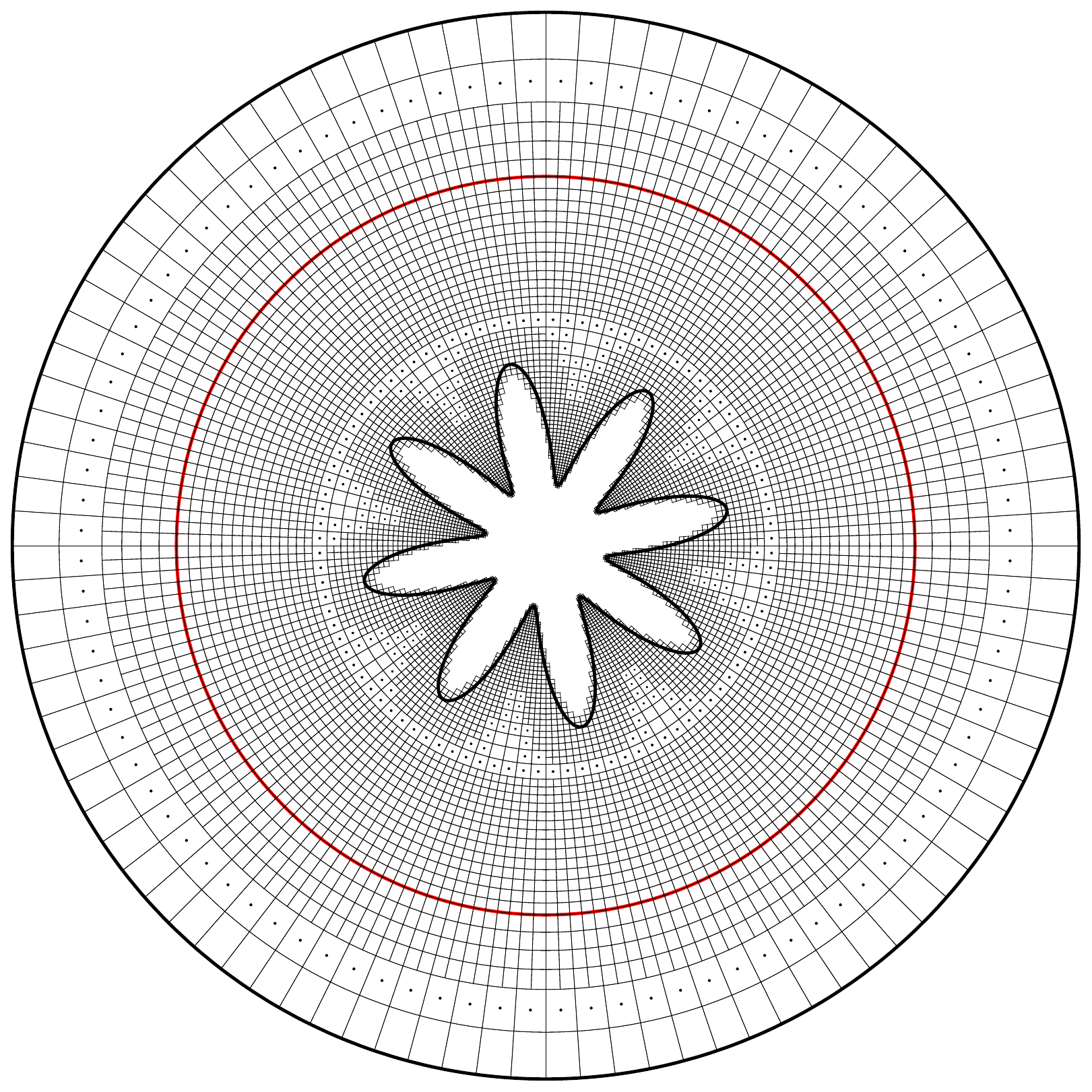}
		\includegraphics[width = 0.3 \linewidth]{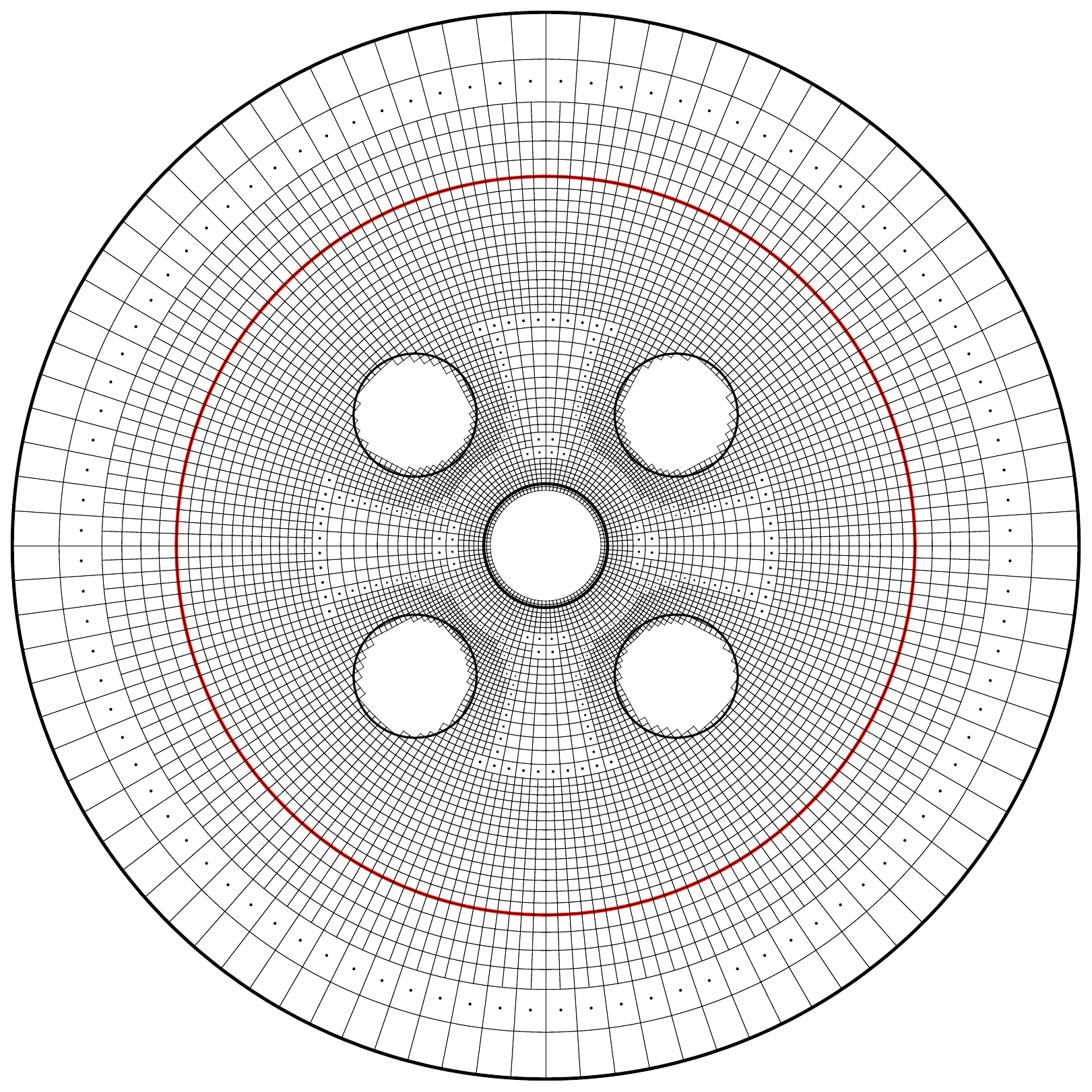}
		\includegraphics[width = 0.3 \linewidth]{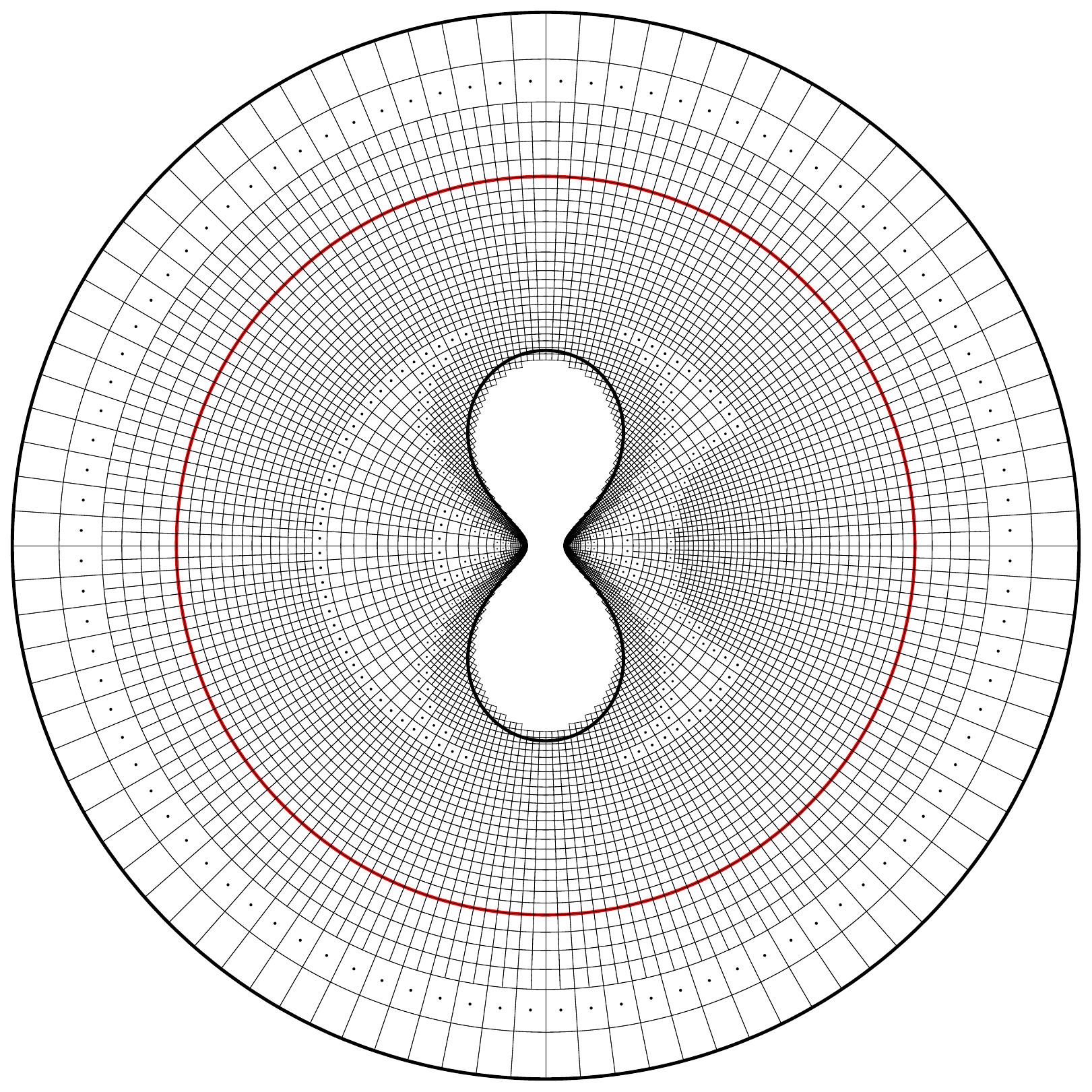}
		
		\vspace{1em}
		\includegraphics[width = 0.31 \linewidth]{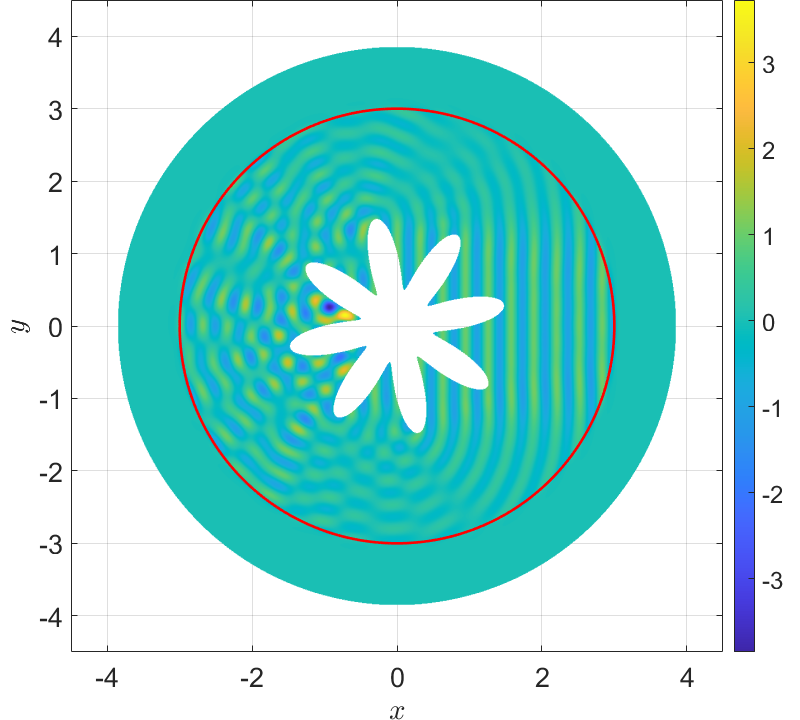}
		\includegraphics[width = 0.31 \linewidth]{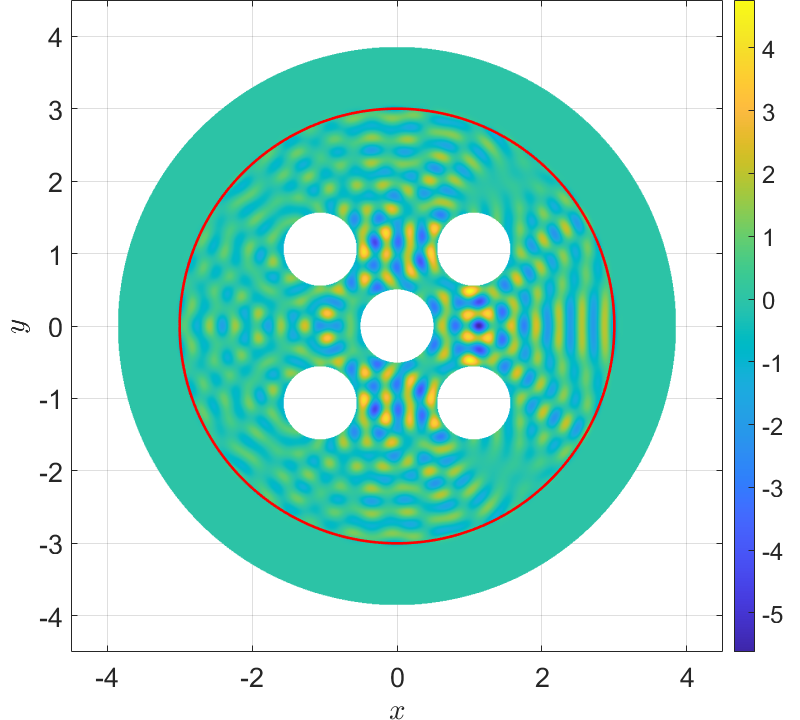}
		\includegraphics[width = 0.31 \linewidth]{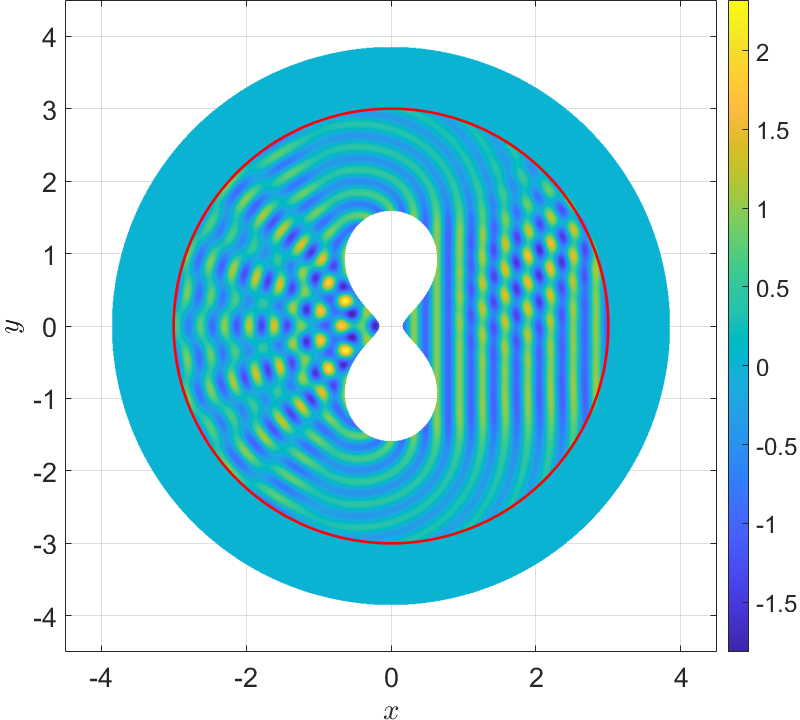}
        \caption{Mesh generated with $\kappa h \approx 1.963$ and real parts of reference solutions in each scenario of \Cref{ex:ex3}. Dots represent the auxiliary nodes and red circle indicates the interface.}
        \label{fig:ex3:mesh_ref_sol}
	\end{figure}

	\newpage
	\begin{figure}[b]
		\vspace{-1em}
	    \centering
		\includegraphics[width = 0.31 \linewidth]{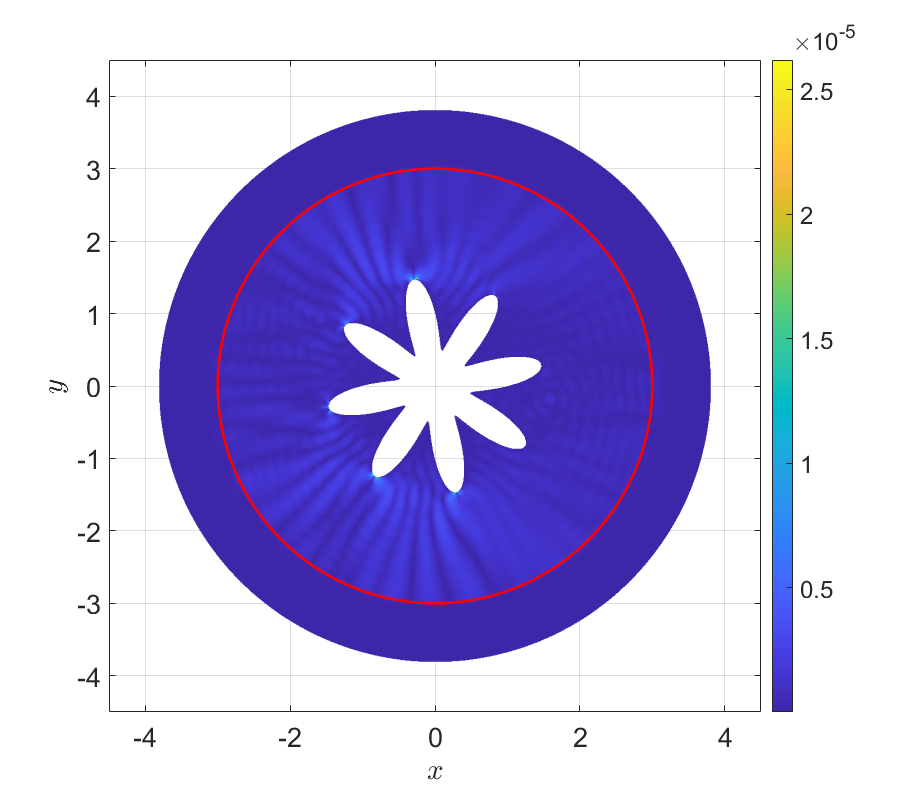}
		\includegraphics[width = 0.31 \linewidth]{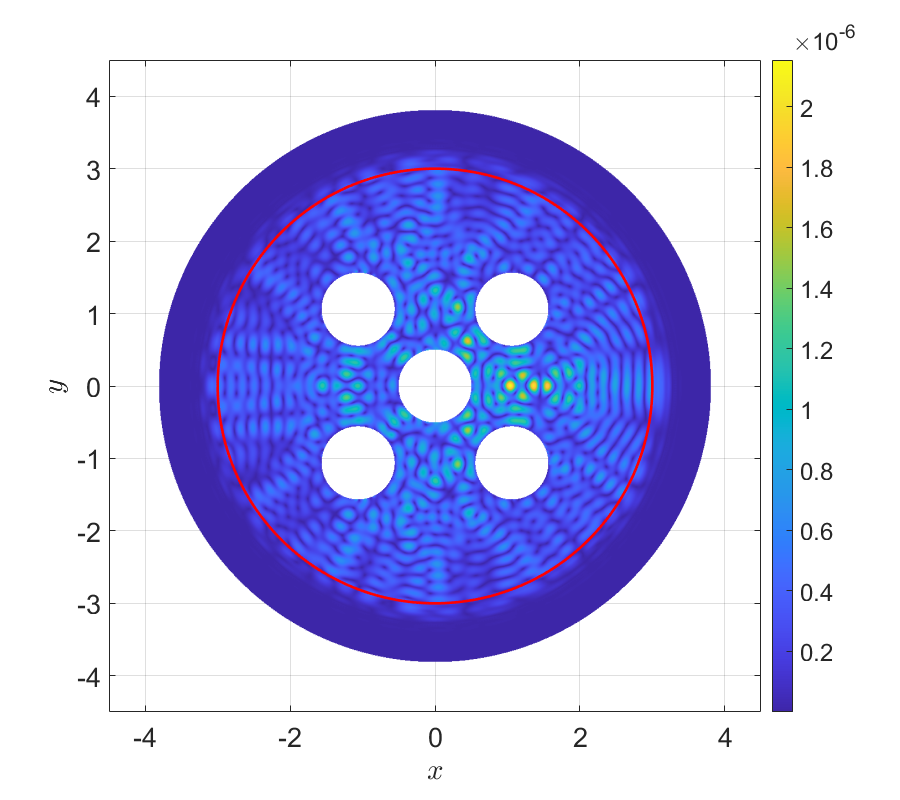}
		\includegraphics[width = 0.31 \linewidth]{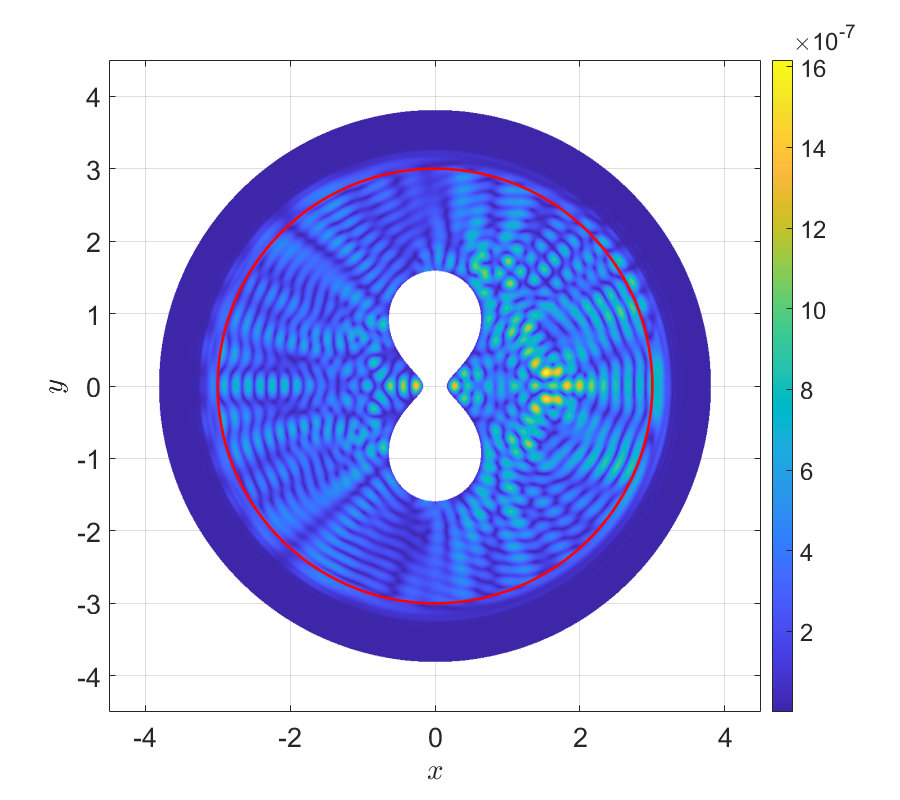}
		
		\vspace{-0.5em}
        \caption{Numerical error $|v_h - v^{\mathrm{ref}}|$ with $\kappa h \approx 0.327$ in \Cref{ex:ex3}.}
        \label{fig:ex3:err}
	\end{figure}
	
	\begin{table}[b]
		\centering
        \begin{NiceTabular}{|c||c|c|c|c|}[cell-space-limits=3pt]
            \hline
            $\kappa h$ & $\ell^2$ error & order & $\ell^\infty$ error & order \\ \hline \hline
            1.963      & 1.013E$-$1     & \phantom{10.00} & 3.821E$-$1          & \phantom{10.00} \\ \hline
            1.309      & 5.215E$-$3     & 7.32  & 4.944E$-$2          & 5.04  \\ \hline
            0.982      & 5.117E$-$4     & 8.07  & 7.001E$-$3          & 6.79  \\ \hline
            0.654      & 5.563E$-$5     & 5.47  & 8.568E$-$4          & 5.18  \\ \hline
            0.491      & 3.978E$-$6     & 9.17  & 1.340E$-$4          & 6.45  \\ \hline
			0.327      & 7.362E$-$7     & 4.16  & 2.590E$-$5          & 4.05  \\ \hline
            0.164      & \Block{1-4}{Reference solution} & & &  \\ \hline \hline
            \Block{1-5}{Without pollution minimization} & & & & \\ \hline
            0.327      & 1.571E$-$4     &  & 2.116E$-$2          &       \\ \hline
        \end{NiceTabular}
        \begin{NiceTabular}{|c|c|c|c|}[cell-space-limits=3pt]
            \hline
            $\ell^2$ error & order & $\ell^\infty$ error & order \\ \hline \hline
            7.176E$-$2     &       & 2.748E$-$1          &       \\ \hline
            6.551E$-$4     & 11.58 & 3.120E$-$3          & 11.04 \\ \hline
            1.078E$-$4     & 6.27  & 5.909E$-$4          & 5.78  \\ \hline
            7.385E$-$6     & 6.61  & 3.728E$-$5          & 6.82  \\ \hline
            1.334E$-$6     & 5.95  & 1.079E$-$5          & 4.31  \\ \hline
			3.894E$-$7     & 3.04  & 2.155E$-$6          & 3.97  \\ \hline
            \Block{1-4}{Reference solution} & & & \\ \hline \hline
            \Block{1-4}{Without pollution minimization} & & &    \\ \hline
            2.126E$-$4     &       & 3.222E$-$3          &       \\ \hline
        \end{NiceTabular}
        
        \begin{NiceTabular}{|c||c|c|c|c|}[cell-space-limits=3pt]
            \hline
            $\kappa h$ & $\ell^2$ error & order & $\ell^\infty$ error & order \\ \hline \hline
            1.963      & 4.465E$-$1     &       & 2.390E$+$0          &       \\ \hline
            1.309      & 3.729E$-$4     & 17.48 & 1.827E$-$3          & 17.70 \\ \hline
            0.982      & 5.765E$-$5     & 6.49  & 2.237E$-$4          & 7.30  \\ \hline
            0.654      & 3.929E$-$6     & 6.62  & 1.396E$-$5          & 6.84  \\ \hline
            0.491      & 6.504E$-$7     & 6.25  & 2.278E$-$6          & 6.30  \\ \hline
			0.327      & 1.527E$-$7     & 3.57  & 1.617E$-$6          & 0.85  \\ \hline
            0.164      & \Block{1-4}{Reference solution} & & &  \\ \hline \hline
            \Block{1-5}{Without pollution minimization} & & & & \\ \hline
            0.327      & 5.169E$-$5     &       & 9.046E$-$4          &       \\ \hline
        \end{NiceTabular}
        \phantom{%
        \begin{NiceTabular}{|c|c|c|c|}[cell-space-limits=3pt]
            \hline
            $\ell^2$ error & order & $\ell^\infty$ error & order \\ \hline \hline
            7.176E$-$2     &       & 2.748E$-$1          &       \\ \hline
            6.551E$-$4     & 11.58 & 3.120E$-$3          & 11.04 \\ \hline
            1.078E$-$4     & 6.27  & 5.909E$-$4          & 5.78  \\ \hline
            7.385E$-$6     & 6.61  & 3.728E$-$5          & 6.82  \\ \hline
            1.334E$-$6     & 5.95  & 1.079E$-$5          & 4.31  \\ \hline
			3.894E$-$7     & 3.04  & 2.155E$-$6          & 3.97  \\ \hline
            \Block{1-4}{Reference solution} & & & \\ \hline \hline
            \Block{1-4}{Without pollution minimization} & & &    \\ \hline
            2.126E$-$4     &       & 3.222E$-$3          &       \\ \hline
        \end{NiceTabular}}
    
		\caption{Numerical error $e_h$ and convergence order in \Cref{ex:ex3}. All data except for the last row are gathered with the pollution minimization process.}
		\label{table:ex3}
    \end{table}
    	
    \begin{example}
    	\label{ex:ex3} \normalfont
    	In this example, we test our FDM with various non-circular scatterers. We keep the same as \Cref{ex:ex1} for $\kappa$, $h^{\mathrm{ref}}$, $h$ (except that we also include $N = 12$), as well as the boundary data $g$. The PML layer is adjusted from $\tilde{r}_* = 3$ and $\tilde{r}_M = 4$, which results in $\kappa d \approx 18.0$ after adjustment. We test three different scenarios: 
		\begin{itemize}
			\item[(a)] $D = \{ (r, \theta): r < 1 + 0.5 \sin(8 \theta) \}$ with $f = 0$;
    		
			\item[(a)] $D = \cup_{(x_0, y_0) \in \mathcal{C}} \{ (x, y): (x - x_0)^2 + (y - y_0)^2 < 0.5^2 \}$, $\mathcal{C} = \{ (0, 0) \} \cup \big\{ (\sin \phi, \cos \phi): \phi = \frac{\pi}{4}, \frac{3\pi}{4}, \frac{5\pi}{4}, \frac{7\pi}{4} \big\}$ with $f = 0$;
    		    		
			\item[(c)] $D = \{ (x^2 + (y - 1)^2) (x^2 + (y + 1)^2) < 0.6 \}$ with source term
	    	\begin{equation*}
				f(r' \cos \theta' + 1.5, r' \sin \theta') = \kappa \exp \big( \ii \kappa r' \cos(\tfrac{\pi}{3} - \theta') \big) \cdot \boldone_{\{ r' < 1 \}} \exp \big( \tfrac{(r')^2}{(r')^2 - 1} \big).
    		\end{equation*}
		\end{itemize}
    	The exponentially stretched mesh with additional mesh refinement near the scatterer $D$ and $\supp f$ is used throughout this example. The mesh corresponding to $N = 12$ (that is, $\kappa h \approx 1.963$) and the reference solutions of each scenario are presented in \Cref{fig:ex3:mesh_ref_sol}. \Cref{table:ex3} shows the numerical errors under $\ell^2$ and $\ell^\infty$ norm and compares the FDMs with and without the pollution minimization process. The $\ell_2$ error is defined by
    	\begin{equation*}
    		e_h := \|v_h - v^{\mathrm{ref}}\|_{\ell^2} = \bigg( \frac{1}{\# \Omega_h} \sum_{\spt \in \Omega_h} |v_h(\spt) - v^{\mathrm{ref}}(\spt)|^2 \bigg)^{1/2}.
    	\end{equation*}
    	Finally, we present the errors $|v_h - v^{\mathrm{ref}}|$ corresponding to $N = 2$ (that is, $\kappa h \approx 0.327$) in \Cref{fig:ex3:err}.
    	
    	To account for the sub-optimal convergence orders in \Cref{table:ex3}, recall that the local truncation error is capped below by a magnitude of $10^{-7}$. Furthermore, we require the mesh to be uniform in $(s, \theta)$ coordinates near the scatterer, which results in a nonuniform mesh in the Cartesian coordinates. Hence, as $h$ gets smaller, there will be a gradually increasing portion of nodes where the threshold $10^{-7}$ of the local truncation error is reached. This causes the deviation from the theoretical convergence order. Despite this, it is true that pollution minimization is still able to improve the accuracy of the results by multiple magnitudes.
    \end{example}

	\appendix
	
	\section{Mathematical justification}
	\label{app:math}
	
	\subsection{Estimates on the solution}
	\label{sec:estimates}
	
	In this section, we establish the decaying property for the exact solution $u$ of \cref{eq:PDE_r} as well as its derivatives. This is beneficial for us to understand the development of our FDM. We will frequently use the Bessel functions $J_j$ and the Hankel functions of the first kind $H_j^{(1)}$, $j \in \Z$. For the definition and properties of these functions, please refer to~\cite{olver2010nist}.

    \begin{lemma}
        \label{lem:H1_H0}
        For any $z \in \C \bs \{0\}$ we have
        \begin{equation}
            \label{eq:H1_H0}
            |H^{(1)}_1 (z)| \leq \left( 1 + C |z|^{-1} \right) |H^{(1)}_0 (z)|,
        \end{equation}
        where $C$ is an absolute constant.
    \end{lemma}

    \begin{proof}
        Note that $H^{(1)}_\nu$ is smooth and has no zeros in the complex plane. Thus we only need to consider when $z \to 0$ and $z \to \infty$. From the asymptotic expansion at 0 \cite[equation (10.7.2), (10.7.7)]{olver2010nist} we obtain
        \begin{equation*}
            |H^{(1)}_1 (z)| \leq C |z|^{-1} |\log z|^{-1} |H^{(1)}_0 (z)|
        \end{equation*}
        for sufficiently small $z$. This obviously implies \eqref{eq:H1_H0} when $z$ is small. Meanwhile, the asymptotic expansion at $\infty$ with error estimates \cite[equation (10.17.13)]{olver2010nist} yields
        \begin{equation*}
            \frac{|H^{(1)}_1 (z)|}{|H^{(1)}_0 (z)|} - 1
            \leq \frac{\tfrac{1}{2} |z|^{-1} + |R_2^+ (0, z)| + |R_2^+ (1, z)|}{|1 - |R_1^+ (0, z)||},
        \end{equation*}
        where $|R_k^+ (\nu, z)| \leq C |z|^{-k}$ for $k = 1$, $2$ and $\nu = 0$, $1$ with an absolute constant $C$. This proves \eqref{eq:H1_H0} when $|z|$ is sufficiently large. 
    \end{proof}

    \begin{lemma}
        \label{lem:Hankel_derivs}
        Let $r_0 > 0$, $j \in \Z$ and $m \in \N_0$. For any $z \in \C$ such that $|z| \geq r_0$, $\re z \geq 0$ and $\im z \geq 0$ we have
        \begin{equation}
            \label{eq:Hankel_derivs}
            \left| \dd{z}[m] H^{(1)}_j(z) \right|
            \leq C_m \exp \! \left( -\kappa \im z \left( 1 - r_0^2 / |z|^2 \right)^{1/2} \right) \left( 1 + \frac{|j|^m + 1}{|z|^m} \right) |H^{(1)}_j (r_0)|.
        \end{equation}
    \end{lemma}

    \begin{proof}
        Since $H^{(1)}_{-j} (z) = e^{\ii j \pi} H^{(1)}_j (z)$, we only need to consider $j \in \N_0$. First, we assume $j \geq m$. Using the identity $H^{(1) \prime}_j (z) = H^{(1)}_{j - 1} (z) - \frac{j}{z} H^{(1)}_j (z)$ repeatedly, we can express $\dd{z}[m] H^{(1)}_j(z)$ as a linear combination of $\{ H^{(1)}_{j - m + k}(z): 0 \leq k \leq m \}$. Moreover, the coefficients of $H^{(1)}_{j - m + k}(z)$ are in the form of $P_k(j) / z^k$, where $P_k$ is a polynomial of degree at most $k$. For the case of $1 \leq j < m$, we additionally use the identity $H^{(1) \prime}_0 (z) = -H^{(1)}_1(z)$. Now $\dd{z}[m] H^{(1)}_j(z)$ becomes a linear combination of $\{ H^{(1)}_k(z): 0 \leq k \leq j \}$, whose coefficients are sum of $P_k(j) / z^k$, $0 \leq k \leq m$. By Young's inequality, we conclude that
        \begin{equation*}
            \left| \dd{z}[m] H^{(1)}_j(z) \right|
            \leq C_m \left( 1 + \frac{|j|^m + 1}{|z|^m} \right) \sup_{0 \leq k \leq j} |H^{(1)}_k (z)|, \ \forall j \geq 1, m \in \N_0.
        \end{equation*}
        Now we apply the decaying property of Hankel functions \cite[Lemma 2.2]{chen2005adaptive} and the monotonicity $|H^{(1)}_k (r_0)| \leq |H^{(1)}_j (r_0)|$ for $0 \leq k \leq j$, $r_0 > 0$ (same proof as \cite[equation (2.18)]{chen2005adaptive}) to show that \eqref{eq:Hankel_derivs} holds. For the remaining case $j = 0$, we have $\dd{z}[m] H^{(1)}_0 (z) = -\dd{z}[m - 1] H^{(1)}_1 (z)$. \Cref{eq:Hankel_derivs} follows by replacing $m$ with $m - 1$ in \eqref{eq:Hankel_derivs} and applying \Cref{lem:H1_H0}.
    \end{proof}

    \begin{prop}
        \label{prop:sol_decay} 
        Suppose $u|_\Gamma \in H^m(\T)$ for some $m \in \N_0$. Then for $\bm{k} = (k_1, k_2) \in \N_0^2$ such that $|\bm{k}| \leq m$, the exact solution $u$ to the Helmholtz equation \eqref{eq:PDE_r} satisfies
        \begin{equation}
            \label{eq:sol_decay}
            \| \partial^{\bm{k}} u(z, \cdot) \|_{L^2 (\T)} 
            \leq C \exp \! \left( -\kappa \im z \left( 1 - r_*^2 / |z|^2 \right)^{1/2} \right) (\kappa^{k_1} + |z|^{-k_1}) \|u|_\Gamma\|_{H^{|\bm{k}|} (\T)}.
        \end{equation}
    \end{prop}

    \begin{proof}
        We write $\bm{k} = (k_1, k_2)$ as usual and start from the series expansion \eqref{eq:u_series}. Taking termwise derivatives and using \Cref{lem:Hankel_derivs}, we obtain
        \begin{align*}
            \| \partial^{\bm{k}} u(z, \theta) \|_{L^2 (\T)}
            &= \left\| \frac{a_n \kappa^{k_1} j^{k_2}}{H^{(1)}_j (\kappa r_*)} \dd{z}[k_1] H^{(1)}_j (\kappa z) \right\|_{\ell^2 (\Z)} \\
            & \leq C_k \exp \! \left( -\kappa \im z \left( 1 - r_*^2 / |z|^2 \right)^{1/2} \right) \left\| a_n j^{k_2} \left( \kappa^{k_1} + \frac{|j|^{k_1} + 1}{|z|^{k_1}} \right) \right\|_{\ell^2 (\Z)} \\
            &\leq C_k \exp \! \left( -\kappa \im z \left( 1 - r_*^2 / |z|^2 \right)^{1/2} \right) (\kappa^{k_1} + |z|^{-k_1}) \|u|_\Gamma\|_{H^{|\bm{k}|} (\T)}.
        \end{align*}
        This also verifies the validity of termwise differentiation.
    \end{proof}
    
	\subsection{Miscellaneous proofs}
	\label{sec:proof}
	
	\begin{lemma}
		\label{lem:center_nz}
		Suppose $\gamma > 0$ and $\rho \in C^1(0, r_M)$ satisfies the assumption ($A_\rho$), then $20 \rho(r)^{-2} + 20 \gamma^2 \rho'(r)^{-2} \neq 0$ for $r \in (0, r_M)$.
	\end{lemma}
	
	\begin{proof}
		We equivalently prove that $\rho'(r) \neq \pm \ii \gamma \rho(r)$ for $r \in (0, r_M)$. This holds for $r_0 \in (0, r_*]$ as $\rho(r)$, $\rho'(r) \in \R$. When $r \in (r_*, r_M)$, we have $\rho'(r) \neq \ii \gamma \rho(r)$ from $\im \rho(r) > 0$ and $\re \rho'(r) > 0$, and $\rho'(r) \neq -\ii \gamma \rho(r)$ from $\re \rho(r) > 0$ and $\im \rho'(r) > 0$.
	\end{proof}
	
	\begin{lemma}
		\label{lem:Bessel_nz}
		Denote $\mathcal{F}_j$ to be either $J_j$ or $H^{(1)}_j$, $j \in \Z$. If $z \in \C \bs \{0\}$ is a common zero of $\mathcal{F}_{j_1}$ and $\mathcal{F}_{j_2}$ for some $j_1, j_2 \in \Z$, then either $|j_1| = |j_2| = 1$ or $|j_1 - j_2| > 2$. In particular, let $\Lambda_z := \{ j \in \Z: \mathcal{F}_j(z) \neq 0 \}$, then for any $z$, $z' \in \C \bs \{0\}$, $\Lambda_z$ has accumulation points at $\pm \infty$ and $\Lambda_z \cap (\Lambda_{z'} + 1)$ is an infinite set.
	\end{lemma}
	
	\begin{proof}
		Suppose there exists $j \in \Z \bs \{0\}$ such that at least two of $\mathcal{F}_{j - 1} (z)$, $\mathcal{F}_j (z)$, $\mathcal{F}_{j + 1} (z)$ are equal to $0$. Then by \cite[equation (10.6.1)]{olver2010nist}, all three of them are $0$. Now \cite[equation (10.6.2)]{olver2010nist} implies $\mathcal{F}_j'(z) = 0$, but it is impossible that both $\mathcal{F}_j(z)$ and $\mathcal{F}_j'(z)$ are $0$ as $\mathcal{F}_j$ satisfies a second order ODE. The same argument holds when two consecutive numbers of $\mathcal{F}_{-1} (z)$, $\mathcal{F}_0 (z)$ and $\mathcal{F}_1 (z)$ are $0$. Hence, we prove the first part of the lemma by contradiction. In other words, we have proved that among any three consecutive integers, at least two of them (except possibly for $\pm 1$) belong to $\Lambda_z$. The second part of this lemma follows naturally.
	\end{proof}
	
	\begin{lemma}
		\label{lem:Bessel_indpt}
		Let $\mathcal{F}_j$ and $\Lambda_z$ be the same as in \Cref{lem:Bessel_nz}, then for any $M \in \N$ and $z$, $z' \in \C \bs \{0\}$, the elements in $E := \big\{ \big( \partial^{\ell_1} \mathcal{F}_j (z') j^{\ell_2} \big)_{j \in \Lambda_z} \in \C^{\#\Lambda_z}: \bm{\ell} \in \LN_{M + 1} \big\}$ are linearly independent. Here $\LN_{M + 1}$ is defined in \Cref{sec:FDM_generic}.
	\end{lemma}
	
	\begin{proof}
		We prove the lemma by contradiction and suppose the elements in $E$ are linearly dependent. Due to \cite[equation (10.6.2)]{olver2010nist} and $\ell_1 \in \{0, 1\}$, we deduce that
		\begin{equation*}
			0 \in {\Span}_{\C} \big\{ \big( \mathcal{F}_{j - \ell_1} (z') j^{\ell_2} \big)_{j \in \Lambda_z}: \bm{\ell} \in \LN_{M + 1} \big\}
			\cap {\Span}_{\C} \big\{ \big( \mathcal{F}_{j + \ell_1} (z') j^{\ell_2} \big)_{j \in \Lambda_z}: \bm{\ell} \in \LN_{M + 1} \big\}.
		\end{equation*}
		In other words, there exists polynomials $P$, $\tilde{P}$, $Q$ and $\tilde{Q}$ with $\deg P$, $\deg \tilde{P} \leq M + 1$, $\deg Q$, $\deg \tilde{Q} \leq M$ such that
		\begin{equation}
			\label{eq:Bessel_recursive}
			P(j) \mathcal{F}_j (z') = Q(j) \mathcal{F}_{j - 1} (z')
			\ \ \text{and} \ \
			\tilde{P}(j) \mathcal{F}_j (z') = \tilde{Q}(j) \mathcal{F}_{j + 1} (z'),
			\ \ \forall j \in \Lambda_z.
		\end{equation}
		
		We now prove a few non-degenerate properties. First, the linear dependence of $E$ implies $P$ and $Q$, or $\tilde{P}$ and $\tilde{Q}$, are not simultaneously zero polynomials. Moreover, if $P$ is a zero polynomial, then $\mathcal{F}_{j - 1} (z') = 0$ for all $j \in \Lambda_z$ except for finitely many zeros of $Q$. In other words, $(\Lambda_z - 1) \cap \Lambda_{z'}$ is a finite set, which contradicts \Cref{lem:Bessel_nz}. Similarly, none of the polynomials in \eqref{eq:Bessel_recursive} is a zero polynomial. Next, suppose there exists $j \in \Lambda_z$ such that any element of $\{ \mathcal{F}_{j'} (z'): j' = j - 1, j, j + 1 \}$ is 0. Using \eqref{eq:Bessel_recursive}, all three elements of the above set are zero, and this again contradicts \Cref{lem:Bessel_nz}.
		
		The above results as well as \Cref{lem:Bessel_nz} show that there is an increasing sequence $\{j_n\}_{n \in \Z}$ such that $j_n \to \pm \infty$ as $n \to \pm \infty$, and $\mathcal{F}_{j_n} (z') / \mathcal{F}_{j_n - 1} (z') = \bo(|j_n|^m)$ for some $m \in \Z$ as $n \to \pm \infty$. On the other hand, the asymptotic expansions \cite[(10.19.1), (10.19.2)]{olver2010nist} and reflection formulas \cite[Section 10.4]{olver2010nist} imply $J_j (z) / J_{j - \mathrm{sign} (j)} (z) = \bo(|j|^{-1})$ and $H^{(1)}_j (z) / H^{(1)}_{j - \mathrm{sign} (j)} (z) = \bo(|j|)$ for $z \in \C \bs \{0\}$ as $j \to \pm \infty$, which is a contradiction. We hence prove the linear independence.
	\end{proof}
	
	\begin{proof}[Proof of \Cref{thm:pollution}]
		Fix $\spt \in \Omega_h = (r, \theta)$. For $p$, $q \in \SS$, define
		\begin{equation*}
			w_{p, q} = \langle w(\spt + ph; \cdot), w(\spt + qh; \cdot) \rangle_{L^2(\T)} + \delta_h \td(p - q).
		\end{equation*}
		Define two matrices $W_h = (w_{p, q})_{p, q \in \SS \bs \{(0, 0)\}}$, $W_h^+ = (w_{p, q})_{p, q \in \SS}$ and a vector $\vec{b}_h = -(w_{p, (0, 0)})_{p \in \SS \bs \{(0, 0)\}}$, then $\mathcal{I}_h (\vec{a}) = \vec{a}^H W_h^+ \vec{a}$, where $\vec{a}^H$ is the conjugate transpose of $\vec{a}$. From either $\delta_h > 0$ in assumption (v$_1$) or assumptions (iii) and (v$_2$), we know that $W_h^+$ is positive definite, so $W_h^+$ as well as $W_h$ are invertible. Besides, the solution for the minimization problem \eqref{eq:hat_C_p} is given by $\widehat{C}_{(0, 0)} (\spt) = 1$ and $(\widehat{C}_p (\spt))_{p \in \SS \bs \{(0, 0)\}} = W_h^{-1} \vec{b}_h$. By the analyticity of $\rho$ in assumption (i) and of $\delta_h$ in assumption (v$_1$) or (v$_2$), $w(\cdot; \theta_0)$ is analytic for all $\theta_0 \in \T$, which implies $w_{p, q}$ is analytic in $h > 0$ for $p, q \in \SS$. Hence, according to \cite[Corollary I]{ribarivc1969analytic}, $(\widehat{C}_p (\spt))_{p \in \SS}$ can be expanded into a Laurent series with respect to~$h$.
		
		Denote $\vec{\tilde{a}}_h := (\widehat{C}_p (\spt))_{p \in \SS}$ and normalize $\vec{\tilde{a}}_h$ into $\vec{a}_h = h^n \vec{\tilde{a}}_h$ with smallest possible $n \in \N$ such that $|\vec{a}_h|_{\ell^\infty} = \bo(1)$, then we have $|\vec{a}_h|_{\ell^\infty} \neq \bo(h)$. For $p \in \SS$ we let $a_p$ be the corresponding component in $\vec{a}_h$ and set $a_p = \sum_{j = 0}^{M + 1} a_{p, j} + \bo(h^{M + 2})$. Now we use \cref{eq:v_taylor_int} and change the order of summation a few times to obtain
		\begin{equation*}
			h^n \widehat{\mathcal{L}}_h w(\cdot; \theta_0) = 
			\sum_{p \in \SS} a_p w(\spt + ph; \theta_0) 
			= \sum_{j = 0}^{M + 1} I_j (\theta_0) h^j + \bo(h^{M + 2})
		\end{equation*}
		for $\theta_0 \in \T$, where
		\begin{equation}
			\label{eq:I_j}
			I_j (\theta_0) := \sum_{\bm{\ell} \in \LN_j} \sum_{k = |\bm{\ell}|}^j \sum_{p \in \SS} a_{p, j - k} A^k_{\bm{\ell}}(p) \partial^{\bm{\ell}} w(\spt; \theta_0).
		\end{equation}
		On the other hand, there exists an $M$-th order consistent FDM with $(C_p(\spt))_{p \in \SS}$ as stencil coefficients by assumption (ii). It follows that
		\begin{equation}
			\label{eq:I_h:comparison}
			\mathcal{I}_h (\vec{a}_h) 
			= h^{2n} \mathcal{I}_h (\vec{\tilde{a}}_h) 
			\leq h^{2n} \mathcal{I}_h \big( (C_p(\spt))_{p \in \SS} \big) 
			= \bo(h^{2M + 4}),
		\end{equation}
		and the last equality comes from the argument at the beginning of \Cref{sec:pollution}. Since $\mathcal{I}_h (\vec{a}_h) = \frac{1}{2\pi} \int_\T |h^n \widehat{\mathcal{L}}_h w(\cdot; \theta_0)|^2 \mathrm{d} \theta_0 + \delta_h |\vec{a}_h|_{\ell^2}^2$, we obtain $h^n \widehat{\mathcal{L}}_h w(\cdot; \theta_0) = \bo(h^{M + 2})$ for $\theta_0 \in \T$. Therefore, $I_j (\theta_0) = 0$ for all $0 \leq j \leq M + 1$ and $\theta_0 \in \T$.
		
		Taking the inverse Fourier coefficients with respect to $\theta_0$ in \eqref{eq:I_j} and eliminating nonzero factors including $\rho'(r)$, we obtain
		\begin{equation}
			\label{eq:I_j_Fourier}
			\sum_{\bm{\ell} = (\ell_1, \ell_2) \in \LN_j} \sum_{k = |\bm{\ell}|}^j \sum_{p \in \SS} a_{p, j - k} A^k_{\bm{\ell}}(p) \cdot \partial^{\ell_1} \mathcal{F}_{j'} (\kappa \rho(r)) (j')^{\ell_2} = 0, \ \forall j \leq M + 1, j' \in \Lambda_*,
		\end{equation}
		where $\Lambda_* := \{ j \in \Z: J_j (\kappa r_*) \neq 0 \}$ and $\mathcal{F}_j = J_j$ when $r \leq r_*$ and $\mathcal{F}_j = H^{(1)}_j$ for $r > r_*$. In particular, this implies
		\begin{equation*}
			0 \in {\Span}_{\C} \big\{ \big( \partial^{\ell_1} \mathcal{F}_j (\kappa \rho(r)) j^{\ell_2} \big)_{j \in \Lambda_*} \in \C^{\#\Lambda_*}: \bm{\ell} \in \LN_{M + 1} \big\}.
		\end{equation*}
		From the linear independence in \Cref{lem:Bessel_indpt}, we return to \cref{eq:I_j_Fourier} and obtain $\sum_{k = |\bm{\ell}|}^j \sum_{p \in \SS}$ $a_{p, j - k} A^k_{\bm{\ell}}(p) = 0$, $\forall j \leq M + 1$, $\bm{\ell} \in \LN_j$. This is equivalent to \eqref{eq:c_p,j} with $c_{p, j}$ replaced by $a_{p, j}$. Since \eqref{eq:c_p,j} is further equivalent to \eqref{eq:FDM_error}, $\vec{a}_h = (a_p)_{p \in \SS}$ also satisfies \eqref{eq:FDM_error}.
		
		It suffices to prove that $(\widehat{C}_p (\spt))_{p \in \SS} = \vec{\tilde{a}}_h$ is the same as $\vec{a}_h$, that is, $n = 0$ when we normalize $\vec{\tilde{a}}_h$ into $\vec{a}_h$. Recall that assumption (ii) implies $\mathcal{I}_h (\vec{\tilde{a}}_h) = \bo(h^{2M + 4})$. In particular, $\delta_h |\vec{\tilde{a}}_h|_{\ell^2}^2 = \bo(h^{2M + 4})$. If assumption (v$_1$) holds, then $\vec{\tilde{a}}_h = \bo(1)$ and it is thus the same as $\vec{a}_h$. Now we assume that (iv) holds. From the expansion $a_p = \sum_{j = 0}^{M + 1} a_{p, j} + \bo(h^{M + 2})$ and the normalization condition in \eqref{eq:hat_C_p} we know that $a_{(0, 0), 0} + \bo(h) = a_{(0, 0)} = h^n \widehat{C}_{(0, 0)} (\spt) = h^n$. On the other hand, $\vec{a}_h$ satisfies \eqref{eq:c_p,j}. By assumption (iv) and $|\vec{a}_h| \neq \bo(h)$, we have $a_{(0, 0), 0} \neq 0$. Hence, we must have $n = 0$. This implies that $(\widehat{C}_p (\spt))_{p \in \SS} = \vec{a}_h$ and it satisfies \eqref{eq:FDM_error}. This finishes the proof.
	\end{proof}
	
	\subsection{Effect of floating point error on pollution minimization process}
	\label{app:float_error}
	
	We aim to deduce the following result informally. We only discuss non-boundary stencils; for boundary stencils we simply replace $M$ and $F(p)$ by $M - 2$ and $F_{\C}(p)$, respectively (see \Cref{sec:non_circular}). 
	
	\vspace{0.5em}\noindent\textit{Claim: Under the notations and assumptions (i), (ii) and (v$_1$) of \Cref{thm:pollution}, if $\delta_h^{-1} = \bo(h^{-2M - 4})$ instead of $\delta_h \sim h^{2M + 4}$, then we can bound the local truncation error as follows:
	\begin{equation}
		\label{eq:float_error}
		\bigg| \mathcal{L}_h v (\spt) - \sum_{p \in \SS} C_p (\spt) F(p) \bigg| \leq C \left( \delta_h^{1/2} + \epsilon_0 \delta_h^{-1/2} \right),
	\end{equation}
	where $\epsilon_0 \approx 10^{-16}$ represents the floating point error.} \vspace{0.5em}
	
	If we do not consider the floating point error, then we can mostly follow the proof of \Cref{thm:pollution} to deduce that the right hand side of \eqref{eq:float_error} is at most $C \delta_h^{1/2}$. The only difference is $\mathcal{I}_h (\vec{a}_h) = \bo(\delta_h)$ in \eqref{eq:I_h:comparison}, which implies that $I_j (\theta_0) = 0$ for all $0 \leq j \leq M' + 1$ and $\theta_0 \in \T$, where $M'$ is the smallest integer such that $\delta_h = \bo(h^{2M' + 4})$. The proof continues by changing $M$ with $M'$.
	
	Now we take the floating point error into consideration. The floating point error may come from multiple sources, but we consider the process of solving the linear system $W_h \vec{a}_h = \vec{b}_h$ as its main source (see the beginning of the proof of \Cref{thm:pollution} for notations).  Due to \eqref{eq:FDM_error}, we have $\min_{|\vec{a}|_{\ell^2} = 1} \widetilde{\mathcal{I}}_h (\vec{a}) \leq C h^{2M + 4}$, where $\widetilde{\mathcal{I}}_h$ is defined in \cref{eq:I_h}. It follows that $\min_{|\vec{a}|_{\ell^2} = 1} \mathcal{I}_h (\vec{a}) \sim \delta_h$. Since $\mathcal{I}_h (\vec{a}) = \vec{a}^H W_h^+ \vec{a}$ and $W_h$ is a principal submatrix of an Hermitian matrix $W_h^+$, we obtain $\| (W_h)^{-1} \| \leq C \delta_h^{-1/2}$. The matrix norm can be taken as the operator norm in $\ell^\infty$ since we are dealing with vectors of bounded length (determined by the stencil size). Hence, the error of the stencil coefficients due to floating point error has an amount of $C \epsilon_0 \delta_h^{-1/2}$. We thus believe that the claim \eqref{eq:float_error} is true.
	
	In view of \eqref{eq:float_error}, if we set $\delta_h \approx \epsilon_0$ when $h$ is sufficiently small, then we achieve the optimal local truncation error of $C \epsilon_0^{1/2}$. In other words, the local truncation error is bounded by $C \max \{ h^{M + 2}, \epsilon_0^{1/2} \}$ for all $h > 0$. In practice, we set $\delta_h$ slightly larger than $\epsilon_0$ (e.g., $10^{-14}$) to mitigate the floating point error when we add the regularization term $\delta_h I$ to the matrix $W_h^+$. 
	
	Beware that the local truncation error is $2$ orders higher than the error of the numerical solution for non-boundary nodes, which means the numerical solution will blow up as $h \to 0$. Therefore, for sufficiently small $h$ (e.g., $\kappa h < 0.15$) we need to return to the generic method in \Cref{sec:FDM_generic}. The pollution minimization process can still be used for any $h > 0$ for boundary stencils, as the local truncation error and the numerical error have the same order. We will consider combining the symbolic computation in the generic method and the pollution minimization process in the future to unify the computation of stencil coefficients. In particular, we may compute the general expression of the stencil coefficients at least for the zeroth and first order terms, and obtain the remainder using a similar minimization process.	

    \bibliographystyle{unsrtnat}
    \bibliography{ref}
    
\end{document}